\documentclass[11pt]{amsart}
\usepackage{amsmath,amscd, xcolor, setspace, amsfonts, bbm, amsthm, tikz, tikz-cd, enumerate, graphicx, mathtools, amssymb}
\usepackage{caption}
\usetikzlibrary{calc,decorations.markings}
\usepackage[margin=1in]{geometry}
\usepackage[shortalphabetic]{amsrefs}
\usepackage{hyperref}
\allowdisplaybreaks
\usepackage{setspace}
\setdisplayskipstretch{2.5}
\newtheorem{thm}{Theorem}[section]
\newtheorem{theorem}[thm]{Theorem}

\newtheorem{lemma}[thm]{Lemma}
\newtheorem*{lemma*}{Lemma}

\newtheorem*{remark*}{Remark}
\numberwithin{equation}{section}
\newcommand{\R}{\mathbb{R}}

\newcommand{\Q}{\mathbb{Q}}
\newcommand{\Z}{\mathbb{Z}}

\newcommand{\C}{\mathbb{C}}

\newcommand{\res}{\underset{\small s=1}{\mathrm{Res}}~}
\newcommand{\sqfree}{~\square\text{-}\mathrm{free}}
\newcommand{\mods}[1]{\,\,(\mathrm{mod}\, ^*{#1})}
\newcommand{\Aut}[1]{|\mathrm{Aut}({#1})|}
\title{A Secondary Term for $D_4$ Quartic Fields Ordered By Conductor}
\author{Matthew Friedrichsen}

\begin{document}
\maketitle

\begin{abstract}
To any quartic $D_4$ extension of $\mathbb{Q}$, one can associate the Artin conductor of a 2-dimensional irreducible representation of the group. Altu\u g, Shankar, Varma, and Wilson determined the asymptotic number of such fields when ordered by conductor. We refine this, realizing a secondary term and power saving error term.
\end{abstract}

\section{Introduction}

For a quartic field $L$ whose normal closure $M/\Q$ has Galois group $D_4,$ Altu\u g, Shankar, Varma, and Wilson \cite{ASVW} define the \textit{conductor} of $L$ (denoted $C(L)$) to be the Artin conductor of the irreducible 2-dimensional Galois representation
$$\rho_M: \text{Gal}(\overline{\Q}/\Q) \to \text{GL}_2(\C)$$
that factors through $\text{Gal}(M/\Q)\cong D_4.$ Using additional algebraic structure coming from an outer automorphism of $D_4$, they proved that if $N_{D_4}( X)$ denotes the number of isomorphism classes of $D_4$ quartic fields with conductor bounded by $X$, then
$$N_{D_4}( X) = \frac{3}{8}\cdot\prod_{p}\left(1 - \frac{1}{p^2} - \frac{2}{p^3} + \frac{2}{p^4}\right) \cdot X \log X + O(X \log \log X).$$

In this paper, we build on their result using an analogue of the Dirichlet hyperbola method that relies on the same algebraic structure. Doing this, we get the following theorem.

\begin{theorem}\label{main_thm}
Let $N_{D_4}( X)$ denote the number of isomorphism classes $[L:\Q]$ where $L$ is a quartic $D_4$ extension with conductor $C(L) \le X$. Then,
\begin{align*}
    N_{D_4}( X) =& \frac{3}{8}\cdot \prod_p \left(1 - \frac{1}{p^2} - \frac{2}{p^3} + \frac{2}{p^4}\right) \cdot X \log X  \\
    &+ \left(\frac{3}{8}\cdot\left(1 - \frac{7\log 2}{20} - 2\sum_{p}\frac{\log p}{p^2 + 2p + 2}\right)\cdot\prod_p \left(1 - \frac{1}{p^2} - \frac{2}{p^3} + \frac{2}{p^4}\right) +c \right)\cdot X\\
    &+ O_{\epsilon}(X^{11/12 + \epsilon})
\end{align*}
where $c$ is given by (\ref{c_defn}).
\end{theorem}

In addition to their main theorem, \cite{ASVW} also proved asymptotics for $D_4$ quartic fields with certain local specifications. They define $\Sigma = (\Sigma_v)_v$ to be a \textit{collection of local specifications} if for each place $v$ of $\Q$, $\Sigma_v$ contains pairs $(L_v,K_v)$ consisting of a degree 4 \' etale alebra $L_v$ of $\Q_v$ and a quadratic subalgebra $K_v$. Like them, we will define the \textit{conductor} to be
$$
    C(L_p, K_p) = \mathrm{Disc}(L_p)/\mathrm{Disc}(K_p).
$$
We will call a collection $\Sigma$ \textit{acceptable} if for all but finitely many places $v$, the set $\Sigma_v$ contains all pairs $(L_v, K_v).$ $\mathcal{L}(\Sigma)$ will denote all $D_4$ quartic fields $L$ such that $L \otimes \Q_v \in \Sigma_v$ for all $v$, and $N_{D_4}(\Sigma; X)$ the number of isomorphism classes in $\mathcal{L}(\Sigma)$ whose conductor is bounded by $X$. To simplify notation we will let
$$\mu(\Sigma_\infty) = \sum_{(L_\infty, K_\infty) \in \Sigma_\infty}\frac{1}{\Aut{L_\infty,K_\infty}}, \text{ and }\mu(\Sigma_p) = \sum_{(L_p, K_p) \in \Sigma_p}\frac{1}{\Aut{L_p,K_p}C(L_p,K_p)},$$
where $\infty$ is the sole infinite place of $\Q$ and $p$ is a prime, and $\Aut{L_v, K_v}$ is the number of automorphisms of $L_v$ that send $K_v$ to itself.

In \cite{ASVW}, quartic $D_4$ fields $L$ are defined to have \textit{central inertia} at odd primes $p$ under certain conditions. This happens locally when $C(L_p, K_p) = p^2$ but $K_p$ is not ramified at $p.$ We will use $\mu(\Sigma_{p^2})$ to denote a sum over pairs $(L_p, K_p)$ having central inertia that is constructed similarly to the sum for $\mu(\Sigma_p).$  As one might imagine, the case for $p=2$ is more complicated. We will investigate this later, but for the moment we will use $\mu(\Sigma_{2^2}), \mu(\Sigma_{2^4}),$ and $\mu(\Sigma_{2^6})$ without explanation. With that, we have the following theorem.

\begin{theorem}\label{local_specs}
If $\Sigma = (\Sigma_v)_v$ is an acceptable collection of local specifications and $m$ is the product of the primes $p$ for which $\Sigma_p$ does not contain every pair $(L_p, K_p)$. When $m \ll X^{1/4},$
\begin{align*}
    N_{D_4}(\Sigma; X) &= \frac{1}{2}X\cdot \left(\log X + 1 - 2\sum_{p}\frac{\log p\cdot\mu(\Sigma_{p^2})}{\mu(\Sigma_p)} - 2\log 2\frac{2 \mu(\Sigma_{2^4}) + 3\mu(\Sigma_{2^6})}{\mu(\Sigma_2)}\right) \cdot \\
    & \quad\mu(\Sigma_\infty)\cdot \prod_p\left(\left(1 - \frac{1}{p}\right)^2\mu(\Sigma_p)\right) + Xc_{\Sigma}  + O_{\epsilon}(X^{11/12 + \epsilon}m^{1/3 + \epsilon}),
\end{align*}
where $c_{\Sigma}$ is a non-multiplicative constant given in (\ref{local_c}). 
\end{theorem}

In order to prove Theorem \ref{local_specs}, we need to count quadratic extensions of some number field $k$ with a set of local specifications for each place of $k$. For a number field $k$ we define $\Sigma_{k} = (\Sigma_{k,v})_v$ to be a \textit{collection of local specifications over} $k$ if for each place $v$ of $k$, $\Sigma_{k,v}$ contains quadratic \'etale alegbras $K_v$ over $k_v.$ Similar to before, we will call $\Sigma_k$ \textit{acceptable} if $\Sigma_{k,v}$ contains all quadratic \'etale algebras $K_v$ except for at possibly finitely many places. Let $\mathcal{K}(\Sigma_k)$ denote all quadratic extensions $K/k$ such that $K \otimes_k k_v \in \Sigma_{k,v}$ and $\Phi_{k,2}(\Sigma_k; C_2, s)$ denote the Dirichlet series
$$\Phi_{k,2}(\Sigma_k; C_2, s) = \sum_{K \in \mathcal{K}(\Sigma_k)} \frac{1}{D_{K/k}^s},$$
where $D_{K/k}$ is the norm of the relative discriminant ideal of $K/k.$

Later on in the paper, we will write down a different formula for this Dirichlet series that will be derived similarly to the method of proving Theorem 1.1 from \cite{CDO}. The formula will give rise to the following theorem.

\begin{theorem}\label{local_asymptotics}
Let $k$ be a number field, $\Sigma_k$ be an acceptable collection of local specifications for $k$. Let $\mathfrak{r}_1$ be the product of odd primes $\mathfrak{p}$ for which $\Sigma_{k,\mathfrak{p}}$ is not complete but contains both $K_{\mathfrak{p}}$ that are ramified at $\mathfrak{p}$ and nothing else, and let $\mathfrak{r}_2$ be the product of odd primes for which $\Sigma_{k,\mathfrak{p}}$ contains only one $K_{\mathfrak{p}}$ ramified at $\mathfrak{p}$ and nothing else. Let $\mathfrak{u}_1$ be the product of odd primes $\mathfrak{p}$ for which $\Sigma_{k, \mathfrak{p}}$ is not complete and contains all $K_{\mathfrak{p}}$ which are not ramified at $\mathfrak{p}$, and let $\mathfrak{u}_2$ be the product of primes for which $\Sigma_{k,2}$ contains one but not both of the $K_{\mathfrak{p}}$ that are unramified at $\mathfrak{p}.$ Then, if $N_{k,C_2}(\Sigma_k; X)$ denotes the number of quadratic extensions $K/k$ such that $K \in \mathcal{K}(\Sigma_k)$ with the norm of the relative discriminant bounded by $X$,
\begin{align*}
    N_{k,C_2}(\Sigma_k; X) =& X\frac{\res \zeta_k(s)}{\zeta_k(2)}\prod_{\sigma_v\mid \infty}\left(\sum_{K_{v} \in \Sigma_{k,v}}\frac{1}{\Aut{K_{v}/k_{v}}}\right)\prod_{\mathfrak{p}}\left(\left(1 + \frac{1}{\mathrm{N}\mathfrak{p}}\right)^{-1}\sum_{K_{\mathfrak{p}} \in \Sigma_{k,\mathfrak{p}}} \frac{1}{\Aut{K_{\mathfrak{p}}/k_{\mathfrak{p}}}D_{K_{\mathfrak{p}}/k_{\mathfrak{p}}}} \right) \\
    & + O_{n,\epsilon}\left(\#\mathrm{Cl}_{\mathfrak{m}}(k)[2]\frac{X^{\frac{n+2}{n+4} + \epsilon}|D_k|^{\frac{1}{n+4} + \epsilon}\mathrm{N}\mathfrak{u}_1^{\epsilon}\mathrm{N}\mathfrak{u}_2^{\frac{1}{n+4}+\epsilon}}{\mathrm{N}\mathfrak{r}_1^{\frac{n+2}{n+4} - \epsilon}\mathrm{N}\mathfrak{r}_2^{\frac{n+1}{n+4} - \epsilon}}\right),
\end{align*}
where $n = [k:\mathbb{Q}], \Aut{K_v/k_v}$ denotes the number of automorphisms of $K_v$ which fix $k_v,$ and $\mathfrak{m} = \mathfrak{m}_\infty \mathfrak{r}_2\mathfrak{u}_2\prod_{\mathfrak{p}\mid 2}\mathfrak{p}^{2e(\mathfrak{p}|2)}$ with $\mathfrak{m}_\infty$ as the product of real infinite embeddings $\sigma_v$ such that $\Sigma_{k,v}$ is not complete at $v$ and $e(\mathfrak{p}|p)$ is the ramification index for $\mathfrak{p}$ over $p.$
\end{theorem}

\subsection{Notation}
For an number field $K$, let $\mathcal{O}_K$ denote its ring of integers and $\mathrm{Cl}(K)$ denote the ideal class group of $K$, $\mathrm{Cl}_{\mathfrak{m}}(K)$ will denote the ray class group with modulus $\mathfrak{m}.$

For any finite extension $L/K,$ we will use $\mathrm{N}_{L/K}\mathfrak{a}$ to denote the ideal norm map $\mathrm{N}_{L/K}:I(L) \to I(K)$ where $I(K)$ is the group of fractional ideals for the number field $K$. We let $\mathrm{N}_{L/K}(\alpha)$ denote the element norm. If $K = \Q$, we may drop the subscript. For ideals, dropping the subscript will mean $\mathrm{N}\mathfrak{a} = |\mathcal{O}_K/\mathfrak{a}|.$ Also, we let $\mathfrak{d}_{L/K}$ denote the relative discriminant ideal in $\mathcal{O}_K$ and $D_{L/K}$ to be the absolute norm of the relative discriminant ideal. We write $D_K$ to denote the absolute discriminant of a number field $K/\Q.$

\section{Setup}
As stated earlier, we will use the Dirichlet hyperbola method to derive Theorem \ref{main_thm}. Recall that for the number of divisors function $\tau(n) = \sum_{a \mid n} 1$, the hyperbola methods gives us
$$ \sum_{n \le X} \tau(n) = \sum_{a \le X^{1/2}}\sum_{b \le X/a} 1 + \sum_{b \le X^{1/2}}\sum_{a \le X/b} 1 - \sum_{a \le X^{1/2}} \sum_{b \le X^{1/2}} 1.$$
Though quartic $D_4$ fields are not as easy to count as the number of divisors on an integer, there are two pieces from \cite{ASVW} that will allow us to use a similar construction. 

The first is Theorem 5.3, which counts the number of quartic $D_4$ fields $L$ up to conductor $X$ with the discriminant of its quadratic subfield $K$ bounded above by $X^{\beta}$ for $0 < \beta < 2/3.$ We will use $\beta = 1/2$ as suggested by the example above. Because the conductor $C(L)$ of a quartic $D_4$ field $L/\Q$ is
$$C(L) = |D_K|\cdot D_{L/K}$$ where $K$ is the quadratic subfield of $L$, we can extend the metaphor from the example above to count these fields. In the first double sum, counting $a$ up to $X^{1/2}$ will instead be counting quadratic fields $K/\Q$ with $|D_K| < X^{1/2}$ and counting $b$ up to $X/a$ will be counting quadratic extensions $L/K$ with $D_{L/K} < X/|D_K|.$ 

However, it's not obvious how to replicate the second double sum from the example and first count $b$ and then $a$. This is where we need the second piece from \cite{ASVW}. Let $M$ denote the Galois closure of $L/\Q$. Then, the second observation is that there is an outer automorphism $\phi \in \text{Gal}(M/\Q)$ such that the fixed field of $\phi(\text{Gal}(M/L))$ is not isomorphic to $L$ but has the same conductor as $L$. They denote this field as $\phi(L)$. Moreover, their Proposition 2.6 implies that if $|D_K| > C(L)^{1/2}$ then if $\phi(K)$ is the quadratic subfield of $\phi(L)$, $|D_{\phi(K)}| < C(L)^{1/2}.$ This means that our equivalent of counting $b \le X^{1/2}$ will be again to count quadratic fields $K/\Q$ with discriminant bounded by $X^{1/2}.$ However, because we are counting isomorphism classes of quartic $D_4$ fields, we will end up multiplying our result by $1/2$ to account for pairs of fields $L/K$ and $L'/K$ that are both quartic $D_4$ fields over $\Q$ and are non-isomorphic over $K$, but are isomorphic over $\Q.$

To complete the metaphor, we need an equivalent of the third double sum. In the original problem, this counts the pairs $(a,b)$ which were counted by both the first and second double sums. However, Proposition 2.6 actually implies a little more. It shows that $D_{L/K} = |D_{\phi(K)}|d^2 2^n$ for some squarefree, odd integer $d$ and some integer $n$. The number $d$ is, in fact, the product of primes $p$ at which $L$ has central inertia. We will also show later on that $n$ must be a non-negative even integer. For convenience, we will say $q=2^i d.$ Thus, in our problem, we have to be careful how we set up this count because the relationship $D_{L/K} = |D_{\phi(K)}|q^2$ complicates which pairs $(K, L)$ get counted multiple times. We will accomplish this by fixing $q < X^{1/4}$ and considering the quartic fields $L/\Q$ such that $D_{L/K}/|D_{\phi(K)}|$ is exactly $q^2$ and then summing over all possible $q$. \cite{ASVW} refers to $D_{L/K}/|D_{\phi(K)}|$ as $J(L)$, which we will also use here. For a specific pair $(K, L)$, this leads us to two possibilities. If $|D_K| < X^{1/2}/q^2$, then the pair will be double counted if $D_{L/K} < X^{1/2}q^2.$ But if $X^{1/2}/q^2 \le |D_K| < X^{1/2}$, then the pair is double counted if $D_{L/K} < X/|D_K|.$ 

With this, we have finished our analogy to the original Dirichlet hyperbola method. So, the sum we want to analyze which counts quartic $D_4$ fields by conductor is
\begin{equation}\label{the_sum}
    \sum_{\substack{[K:\Q]=2\\ |D_K| < X^{1/2}}}\sum_{\substack{[L:K]=2 \\ D_{L/K} < X/|D_K|}} 1 - \frac{1}{2}\sum_{q < X^{1/4}}\left( \sum_{\substack{[K:\Q]=2\\ |D_K| < X^{1/2}/q^2}}\sum_{\substack{[L:K]=2 \\ D_{L/K} < X^{1/2}q^2 \\ J(L)=q^2}} 1 + \sum_{\substack{[K:\Q]=2\\ X^{1/2}/q^2\le|D_K| < X^{1/2}}}\sum_{\substack{[L:K]=2 \\ D_{L/K} < X/|D_K| \\ J(L)=q^2}} 1\right).
\end{equation}

In Section 3, we will prove the relationship $J(L) = q^2$ more generally and discuss the possibilities for the valuation of $q$ at 2. In Section 4, we will prove a theorem similar to Theorem 1.1 from \cite{CDO} which will give us Theorem \ref{local_asymptotics} for counting quadratic extensions $K/k$ with local specifications. Finally, in Section 5, we will prove Theorems \ref{main_thm} and \ref{local_specs}.

\section{The Flipped Field}
Going forward, for a quartic $D_4$ field $L/\mathbb{Q}$ with quadratic subfield $K$, we will call $\phi(K)$ the \textit{flipped field}. If we let, $L_1$ be a $D_4$ field over $\mathbb{Q}$ and $M$ be its normal closure, then in the subfield diagram below, $K_2$ is the flipped field with $\phi(L_1) = L_2.$ $L_1'$ and $L_2'$ are the conjugate fields of $L_1$ and $L_2$, respectively.

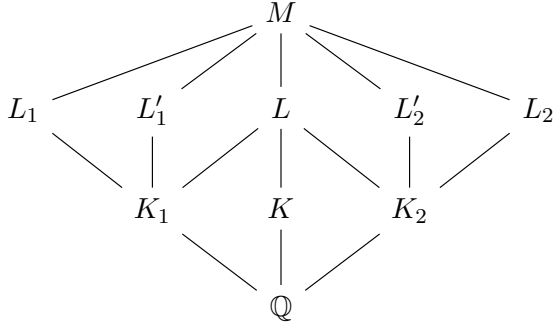
\begin{center}
\begin{tikzcd}
                         &                                                                & M                                                                     &                                                                &                          \\
L_1 \arrow[rru, no head] & L_1' \arrow[ru, no head]                                       & L \arrow[u, no head]                                                  & L_2' \arrow[lu, no head]                                       & L_2 \arrow[llu, no head] \\
                         & K_1 \arrow[ru, no head] \arrow[u, no head] \arrow[lu, no head] & K \arrow[u, no head]                                                  & K_2 \arrow[lu, no head] \arrow[u, no head] \arrow[ru, no head] &                          \\
                         &                                                                & \mathbb{Q} \arrow[lu, no head] \arrow[u, no head] \arrow[ru, no head] &                                                                &                         
\end{tikzcd}
\end{center}

The authors of \cite{ASVW} prove that the $p$-part of $J(L_1) = D_{L_1/K_1}/D_{K_2}$ (or $J_p(L_1)$) is either $p^2$ or $1$ for odd primes $p$ by examining the table of possible decomposition and inertia subgroups of $D_4$ for $p$ and showing that $J_p(L_1) = p^2$ under a certain condition. We will prove this from a different perspective in the lemma below which will also allow us to show that $J_2(L_1)$ is a square.

\begin{lemma}\label{norm_alpha}
Let $k$ be a number field and $L_1/k$ be a quartic $D_4$ extension of $k$. If $K_1$ is the quadratic subfield of $L_1$ and $K_2$ is the flipped field of $L_1$, then there exists an element $\alpha \in K_1$ such that $L_1 = K_1(\sqrt{\alpha})$ and $K_2 = k(\sqrt{\mathrm{N}_{K_1/k}(\alpha)}).$
\end{lemma}

\begin{proof}
Because $L_1/k$ is quartic with Galois closure $D_4$, we know that $L = k(\theta)$ for some $\theta$ with minimal polynomial $f(x) = x^4 + Ax^2 + B \in k[x]$ from Theorem 4.1 of \cite{D4_Galois}. Examining the proof of the theorem, we see that $K_1 = k(\theta^2)$ with $\theta^2 = (-A \pm \sqrt{A^2-4B})/2$ and $A^2 - 4B \ne \square \in k$. Without loss of generality, we can take $\alpha = (-A + \sqrt{A^2 - 4B})/2$ and $L = K(\sqrt{\alpha}).$ We let $\bar{\beta} = (-A - \sqrt{A^2 - 4B})/2.$  So, the complete set of roots for $f(x) = \{ \pm \sqrt{\alpha}, \pm \sqrt{\bar{\beta}}\}.$

Now, from Theorem 1.2 and the proof of Theorem 3.14 in \cite{D4_Galois}, we know that disc $f$ is not a square in $k$ and $K_2 = k(\sqrt{\text{disc}f}).$ But $\text{disc}f = \prod_{i < j} (r_i - r_j)^2$ where $r_i$ are the roots of $f(x).$ For convenience, lets say $r_1 = \sqrt{\alpha},$ $r_2 = -r_1,$ $r_3 = \sqrt{\bar{\beta}},$ and $r_4 = -r_3.$ Then, this product is
\begin{align*}
    \text{disc} f &= 16r_1^2(r_1 - r_3)^4(r_1 + r_3)^4r_3^2 \\
    & = 16r_1^2 r_3^2 \left((r_1- r_3)(r_1 + r_3)\right)^4 \\
    & = 16 \cdot \mathrm{N}_{K_1/k}(\alpha) \cdot (r_1^2 - r_3^2)^4.
\end{align*}
But $(r_1^2 - r_3^2)^4 = (A^2 - 4B)^2,$ which is a square in $k$. Thus, $K_2 \cong k(\sqrt{\mathrm{N}_{K_1/k}(\alpha)}).$
\end{proof}

With this lemma, we can make the more general claim:
\begin{lemma}\label{differ_by_square}
Let $k$ be a number field and $L_1/k$ be a quartic $D_4$ extension of $k$. If $K_1$ is the quadratic subfield of $L_1$ and $K_2$ is the flipped field of $L_1$, then $\mathrm{N}_{K_1/k}\mathfrak{d}_{L_1/K_1} = \mathfrak{d}_{K_2/k}\mathfrak{q}^2$ for some integral ideal $\mathfrak{q}$ of $\mathcal{O}_k.$
\end{lemma}

\begin{proof}
To prove this statement, we will use the work of Cohen, Diaz-y-Diaz, and Olivier \cite{CDO}. Their Proposition 3.4 states that if $L_1 = K_1(\sqrt{\alpha})$ as in Lemma \ref{norm_alpha}, then $\mathfrak{d}_{L_1/K_1} = 4\mathfrak{a}/\mathfrak{c}^2$ where $\mathfrak{a}$ is the largest squarefree ideal dividing $\alpha\mathcal{O}_K$ and $\mathfrak{c}$ is the largest ideal dividing 2 such that $(\mathfrak{a},\mathfrak{c}) = 1$ and $x^2 \equiv \alpha \mods{\mathfrak{c}^2}$ has a solution. Here, $x \equiv y \mods{\mathfrak{a}}$ has the standard meaning from class field theory that $\nu_{\mathfrak{p}}(x-y) \ge \nu_{\mathfrak{p}}(\mathfrak{a})$ for every $\mathfrak{p}$ dividing $\mathfrak{a}.$

For any odd prime $\mathfrak{p}$ of $\mathcal{O}_k$, Proposition 3.4 and Lemma \ref{norm_alpha} together imply $v_{\mathfrak{p}}(\mathrm{N}_{K_1/k}\mathfrak{d}_{L_1/K_1}) \equiv v_{\mathfrak{p}}(\mathfrak{d}_{K_2/k}) \mod 2$ and $v_{\mathfrak{p}}(\mathrm{N}_{K_1/k}\mathfrak{d}_{L_1/K_1}) \ge v_{\mathfrak{p}}(\mathfrak{d}_{K_2/k}).$ So, we need to consider the behavior at primes dividing 2.

As implied by the proof of Proposition 3.4 in \cite{CDO}, for any $\mathfrak{p} \mid 2, v_{\mathfrak{p}}(\mathfrak{d}_{K_2/k})$ is odd if any only if $\mathfrak{p}$ divides the squarefree part of $\mathrm{N}_{K_1/k}(\alpha).$ So, we have $v_{\mathfrak{p}}(\mathrm{N}_{K_1/k}\mathfrak{d}_{L_1/K_1}) \equiv v_{\mathfrak{p}}(\mathfrak{d}_{K_2/k}) \mod 2.$ Moreover, if $v_{\mathfrak{p}}(\mathfrak{d}_{K_2/k})$ is odd, we also have $v_{\mathfrak{p}}(\mathrm{N}_{K_1/k}\mathfrak{d}_{L_1/K_1}) \ge v_{\mathfrak{p}}(\mathfrak{d}_{K_2/k})$ because $v_{\mathfrak{p}}(\mathrm{N}_{K_1/k}\mathfrak{d}_{L_1/K_1}) = 2e(\mathfrak{P}|\mathfrak{p})e(\mathfrak{p}|2) + 1$ for some prime $\mathfrak{P}$ lying over $\mathfrak{p}$ in $\mathcal{O}_{K_1}$ and $v_{\mathfrak{p}}(\mathfrak{d}_{K_2/k}) = 2e(\mathfrak{p}|2) + 1.$

Assume $v_{\mathfrak{p}}(\mathfrak{d}_{K_2/k})$ is even. By Proposition 3.4, we know that $v_{\mathfrak{p}}(\mathfrak{d}_{K_2/k}) = 2(e-m)$ where $e = e(\mathfrak{p}\mid 2)$ and $0 \le m \le e$ such that $\mathrm{N}_{K_1/k}(\alpha)$ is a square mod $\mathfrak{p}^{2m}$ but not mod $\mathfrak{p}^{2(m+1)}$ (except if $m=e$ in which case $\mathfrak{p}$ does not ramify in $K_2$). To show that $v_{\mathfrak{p}}(\mathrm{N}_{K_1/k}\mathfrak{d}_{L_1/K_1}) \ge v_{\mathfrak{p}}(\mathfrak{d}_{K_2/k})$, it will be enough to show that for any $\alpha \in K_1$ and $\mathcal{O}_k$-ideal $\mathfrak{c}$ if $\mathrm{N}_{K_1/k}(\alpha)$ is not a square mod $\mathfrak{c}$, then $\alpha$ is not a square mod $\mathfrak{c}\mathcal{O}_{K_1}.$ We will do this by proving the contrapositive.

Assume $x^2 \equiv \alpha \mods{\mathfrak{c}\mathcal{O}_{K_1}}$, so $\alpha = x^2 + c$ for some $c \in \mathfrak{c}\mathcal{O}_{K_1}$. Now, because $K_1/k$ is quadratic, then $K_1 = k(\sqrt{\beta})$ for some $\beta \in k^{\times}\setminus k^{\times 2}$ and $x = x_1 + x_2\sqrt{\beta}, c = c_1 + c_2\sqrt{\beta}$ for some $x_1, x_2, c_1, c_2 \in k$. Thus
\begin{align*}
    \mathrm{N}_{K_1/k}(\alpha) &= \mathrm{N}_{K_1/k}(x^2 + c) \\
    &= (x_1^2 + x_2^2\beta + c_1)^2 - (2x_1x_2 + c_2)^2\beta \\
    &= \mathrm{N}_{K_1/k}(x^2) + \mathrm{N}_{K_1/k}(c) + 2(x_1^2 + x_2^2\beta)c_1 - 4x_1x_2c_2\beta \\
    &= \mathrm{N}_{K_1/k}(x^2) + \mathrm{N}_{K_1/k}(c) + \mathrm{Tr}_{K_1/k}(\bar{x}^2c)\\
    &\equiv \mathrm{N}_{K_1/k}(x^2) \mods{\mathfrak{c}}.
\end{align*}
For the last two lines $\bar{x} = x_1 - x_2\sqrt{\beta}$ and we have $\mathrm{Tr}_{K_1/k}(\bar{x}^2c) \in \mathfrak{c}$ because $\bar{x}^2c \in \mathfrak{c}\mathcal{O}_{K_1}$ and $K_1/k$ is Galois.
\end{proof}

Now, we can see that when $k= \Q$, the possible values for $J_2(L_1)$ are of the form $2^{2i}$ for $i = 0, 1, 2, 3.$ But, we still don't have a way to count quartic $D_4$ fields where $J_2(L_1) = 2^{2i}.$ For this we will need to think locally.

\subsection{Local Fields and the Flipped Field}
Let $L/\Q$ and $L'/\Q$ be non-isomorphic quartic $D_4$ fields such that $L$ and $L'$ have the same local conditions at 2. Ideally, we would have $J_2(L) = J_2(L')$, but since $J(L) = D_{L/K}/D_{\phi(K)}$, this would require that $\phi(K)$ and $\phi(K')$ have the same local conditions at 2. Fortunately, this is implied by the following lemma.
\begin{lemma}\label{local_flipped}
Let $k$ be a number field, $\mathfrak{p}$ be a prime ideal of $\mathcal{O}_k$, and $L/k$ and $L'/k$ be two non-isomorphic $D_4$ quartic extensions such that $L \otimes_{k} k_{\mathfrak{p}} \cong L' \otimes_{k} k_{\mathfrak{p}}.$ Then, $\phi(K) \otimes_{k} k_{\mathfrak{p}} \cong \phi(K') \otimes_{k} k_{\mathfrak{p}}.$
\end{lemma}

Before proving the lemma, let's consider a simpler case. If we have two quadratic extensions $K = k(\sqrt{\alpha})$ and $K' = k(\sqrt{\alpha'})$ and some prime ideal $\mathfrak{p}$ of $\mathcal{O}_k$, then $K \otimes_{k} k_{\mathfrak{p}} \cong K' \otimes_{k} k_{\mathfrak{p}}$ if and only if $\alpha$ and $\alpha'$ are in the same class of $k_{\mathfrak{p}}^{\times}/k_{\mathfrak{p}}^{\times 2}.$ We will use this idea in the proof.
\begin{proof}
If $K$ and $K'$ are the quadratic subfields of $L$ and $L'$, respectively, then we know that $$K \otimes_k k_{\mathfrak{p}} \cong K' \otimes_k k_{\mathfrak{p}} \cong \prod_{\substack{\mathfrak{P} \subset \mathcal{O}_K\\ \mathfrak{P}\mid \mathfrak{p}}} K_{\mathfrak{P}}.$$ 
Assume $L= K(\sqrt{\beta})$ and $L'=K'(\sqrt{\beta'})$ for some $\beta \in K$ and $\beta' \in K'.$ Let $\Phi_K, \Phi_{K'}$ be the isomorphisms from $K \otimes_{k} k_{\mathfrak{p}}$ and $K' \otimes_{k} k_{\mathfrak{p}}$ to $\prod K_{\mathfrak{P}},$ and $\Phi_{K,\mathfrak{P}}, \Phi_{K', \mathfrak{P}}$ be the restriction maps to $K_{\mathfrak{P}}.$ Then, because $L \otimes_{k} k_{\mathfrak{p}} \cong L' \otimes_{k} k_{\mathfrak{p}}$, we must have for each $\mathfrak{P}$ that $\Phi_{K,\mathfrak{P}}(\beta) = \Phi_{K', \mathfrak{P}}(\beta')a^2$ for some $a \in K_{\mathfrak{P}}.$ So $\Phi_{K}(\beta) = \Phi_{K'}(\beta)'\bar{a}^2$ for some $\bar{a} \in \prod K_{\mathfrak{P}}.$ Thus $N_{K \otimes_{k} k_{\mathfrak{p}}/k_{\mathfrak{p}}}(\beta) = N_{K' \otimes_{k} k_{\mathfrak{p}}/k_{\mathfrak{p}}}(\beta')b^2$ for some $b \in k_{\mathfrak{p}}$.
\end{proof}

With this lemma, we can finish classifying the possibilities for $J_2(L)$. Since $L = K(\sqrt{\alpha})$, we can go further and use Serre's mass formula \cite{Serre} to assign weights to the value of $J_2(L)$ based on the class of $\alpha$ in $K_{\mathfrak{p}}^{\times}/K_{\mathfrak{p}}^{\times 2}$ for each $\mathfrak{p} \mid 2$ in $\mathcal{O}_K.$ If we fix a quadratic \'etale algebra $K_2 = K \otimes_{\Q} \Q_2$ and compute
$$\sum_{\substack{(L_2,K_2)\\J_2(L) = 2^{2i}}}\frac{1}{\Aut{L_2,K_2}C(L_2,K_2)},$$
we get the total weight for all $L/K$ such that $J_2(L) = 2^{2i}$ where $K \otimes_{\Q}\Q_2 \cong K_2.$ We present the weights in the tables below organized by the valuation at 2 of $D_{L/K}$ and $D_{\phi(K)}.$ Though there are 8 different isomorphism classes for degree 2 \'etale algebras over $\Q_2$, we found that these can be reduced to 4 cases for $K/\Q$. The code to compute the tables is at \href{https://github.com/friedrichsenm/d4_by_conductor}{https://github.com/friedrichsenm/d4\_by\_conductor}. Tables 3 and 4 are computed by adding together the results for all cases where $D_K$ matches the congruence condition.

\begin{center}
    \begin{tabular}{|l|c|c|c|}
        \hline
         & $v_2(D_{\phi(K)})=$0 & 2 & 3 \\
         \hline
         $v_2(D_{L/K})=$0&1/2 &0 &0 \\
         2 & 0 &1/4 & 0 \\
         3 & 0& 0 &1/4 \\
         4 & 1/32& 0& 0\\
         5 & 0&0 &1/16 \\
         6 & 1/64&1/64 &0 \\
         \hline
    \end{tabular}
    \captionof{table}{Weights of $v_2(D_{L/K})$ vs $v_2(D_{\phi(K)})$ when $D_K \equiv 1$ mod $8$}
\end{center}

\begin{center}
    \begin{tabular}{|l|c|c|}
        \hline
         & $v_2(D_{\phi(K)})=$0 & 2 \\
         \hline
         $v_2(D_{L/K})=$0&1/2 &0 \\
         4 & 1/32& 1/16\\
         6 & 1/64&1/64 \\
         \hline
    \end{tabular}
    \captionof{table}{Weights of $v_2(D_{L/K})$ vs $v_2(D_{\phi(K)})$ when $D_K \equiv 5$ mod $8$}
\end{center}

\begin{center}
    \begin{tabular}{|l|c|c|}
        \hline
         & $v_2(D_{\phi(K)})=$0  & 3 \\
         \hline
         $v_2(D_{L/K})=$0&1/4 & 0\\
         2 & 1/16&  0\\
         4 & 1/32& 0\\
         5 & 0 &1/32 \\
         \hline
    \end{tabular}
    \captionof{table}{Weights of $v_2(D_{L/K})$ vs $v_2(D_{\phi(K)})$ when $D_K \equiv 4$ mod $8$}
\end{center}

\begin{center}
    \begin{tabular}{|l|c|c|c|}
        \hline
         & $v_2(D_{\phi(K)})=$0 & 2 & 3 \\
         \hline
         $v_2(D_{L/K})=$0&1/4 &0 & 0\\
         2 & 1/16&0 &  0\\
         4 & 0& 1/32& 0\\
         5 & 0& 0 &1/32 \\
         \hline
    \end{tabular}
    \captionof{table}{Weights of $v_2(D_{L/K})$ vs $v_2(D_{\phi(K)})$ when $D_K \equiv 0$ mod $8$}
\end{center}

\section{Counting Quadratic Fields with Local Conditions}
As alluded to earlier on, we need to prove that the weights calculated using the mass formula match with weights when counting global fields $L/K$ that match the local conditions. To do this we will prove a theorem in the vein of Theorem 1.1 from \cite{CDO} that will allow us to single out specific local conditions at 2 we wish to include. Similar to \cite{CDO}, we let
$$\Phi_{k,2}(C_2,s) = \sum_{[K:k]=2} \frac{1}{D_{K/k}^s}.$$ Then, we have the following theorem.
\begin{theorem}\label{local_CDO}
Let $k$ be a number field, $\mathfrak{n}$ be the maximal squarefree ideal divisor of $2$ and $\mathfrak{m} = \prod_{\mathfrak{p} \mid 2} \mathfrak{p}^{2e(\mathfrak{p}|2) +1}.$ For each $\mathfrak{p} \mid 2$, select a uniformizer $\pi_\mathfrak{p} \in \mathcal{O}_k$ such that $\pi_\mathfrak{p} \equiv 1 \mods{\mathfrak{q}}$ for every $\mathfrak{q} \mid 2$ where $\mathfrak{p} \ne \mathfrak{q}.$  Then, the Dirichlet series $\Phi_{k, 2}(C_2, s)$ is
$$\Phi_{k,2}(C_2,s) = -1 + \frac{1}{2^{g + i(k)}\zeta_k(2s)}\prod_{\mathfrak{p} \mid 2}\left(1 - \mathrm{N}\mathfrak{p}^{-2s}\right)^{-1} \sum_{\mathfrak{c}\mid \mathfrak{n}} \sum_{\bar{\beta} \in \mathcal{O}_{k,\mathfrak{m}}/\mathcal{O}_{k,\mathfrak{m}}^2} D_2(\pi_{\mathfrak{c}} \beta)^{-s}\sum_{\chi \in \mathrm{Cl_{\mathfrak{m}}^2(k)}} \chi(\pi_{\mathfrak{c}}^{-1}\beta^{-1}\mathfrak{c}) L_k(s,\chi),$$
where $g$ is the number of prime ideals in $\mathcal{O}_k$ dividing $2$, $\mathcal{O}_{k,\mathfrak{m}} = (\mathcal{O}_k/\mathfrak{m})^{\times}$, $D_2(x) = \prod_{\mathfrak{p} \mid 2}D_{k_{\mathfrak{p}}(\sqrt{x})/k_{\mathfrak{p}}},$ and $\pi_{\mathfrak{c}} = \prod_{\mathfrak{p} \mid \mathfrak{c}} \pi_\mathfrak{p},$ and $\beta$ is a lift of $\bar{\beta}$ to $\mathcal{O}_{k,\mathfrak{m}}.$
\end{theorem}

In order to prove the theorem, we will use many of the same objects as in \cite{CDO}. We will recall some of their definitions here: 
\begin{itemize}
    \item An element $u \in k^{\times}$ is a \textit{virtual unit} if there exists an ideal $\mathfrak{q}$ such that $u\mathcal{O}_k = \mathfrak{q}^2$. It is clear that the set of virtual units $V(k)$ is a group.
    \item We will call the quotient group $S(k) = V(k)/k^{\times 2}$ the \textit{Selmer group} of $k$.
    \item Let $\mathfrak{m} = \mathfrak{m}_0\mathfrak{m}_{\infty}$ be a modulus of $k$. The \textit{ray Selmer group modulo} $\mathfrak{m}$ is the subgroup $S_{\mathfrak{m}}(k)$ of $S(k)$ of elements $\bar{u}$ such that for some lift $u$ of $\bar{u}$ coprime to $\mathfrak{m}_0$ there exists a solution to $x^2 \equiv u \mods{\mathfrak{m}_0}$ and for each real infinite prime $\sigma \mid \mathfrak{m}_{\infty}$ we have $\sigma(u) > 0.$
\end{itemize}
Note that this last definition is different from the one in \cite{CDO} and we will want to reprove some results using this modified definiton. First, a modified version of Lemma 3.5 from \cite{CDO}:
\begin{lemma}\label{ray_class_square}
Let $\mathfrak{m} = \mathfrak{m}_0\mathfrak{m}_{\infty}$ be a modulus of $k$, and let $\mathfrak{a}$ be an integral ideal coprime to $\mathfrak{m}_0$ such that there exists an ideal $\mathfrak{q}$ also coprime to $\mathfrak{m}_0$ with $\mathfrak{aq}^2 = \alpha_0\mathcal{O}_k.$ The following two conditions are equivalent.
\begin{enumerate}
    \item There exists an element $\bar{u}$ of $S(k)$ such that, for any lift $u$ of $\bar{u}$ coprime to $\mathfrak{m}_0$, the congruence $x^2 \equiv \alpha_0u \mods{\mathfrak{m}_0}$ has a solution and $\sigma(\alpha_0u) > 0$ for every $\sigma \mid \mathfrak{m}_{\infty}.$
    \item The class of $\mathfrak{a}$ is a square in the ray class group $\mathrm{Cl}_{\mathfrak{m}}(k).$
\end{enumerate}
\end{lemma}

\begin{proof}
The proof will largely follow the proof of Lemma 3.5 in \cite{CDO} with a few changes. Assume (1). Then $x^2 = \alpha_0u\beta$ with $\beta \equiv 1 \mods{\mathfrak{m}}$ and $\sigma(\beta) > 0$ for every $\sigma \mid \mathfrak{m}_{\infty}.$ The rest of this direction of the proof is the same as in the cited proof.

Now, assume (2). The only part of this direction of the proof that is different is noting that $\sigma(\beta') > 0$ for every $\sigma \mid \mathfrak{m}_{\infty}$ and that because $\alpha_0u = \beta'$ , then $\sigma(\alpha_0u) > 0$ for every $\sigma \mid \mathfrak{m}_{\infty}.$
\end{proof}

Next an equivalent of Lemma 3.7 from the same paper.

\begin{lemma}
Let $\mathfrak{m}$ be a modulus as before and for notational simplicity, let $\mathcal{O}_{k,\mathfrak{m}} = (\mathcal{O}_k/\mathfrak{m}_0)^{\times} \times \{\pm 1\}^{|\mathfrak{m}_{\infty}|}.$ Then, the following sequence is exact.
$$1 \longrightarrow S_{\mathfrak{m}}(k) \longrightarrow S(k) \longrightarrow \mathcal{O}_{k,\mathfrak{m}}/\mathcal{O}_{k,\mathfrak{m}}^2 \longrightarrow \mathrm{Cl}_{\mathfrak{m}}(k)/\mathrm{Cl}_{\mathfrak{m}}(k)^2 \longrightarrow \mathrm{Cl}(k)/\mathrm{Cl}(k)^2 \longrightarrow 1.$$
\end{lemma}

\begin{proof}
The proof for this lemma is the same except to note that the element $\beta$ used in the proof will also be positive for every $\sigma \mid \mathfrak{m}_{\infty}.$
\end{proof}

We will also be concerned with the cardinality of $S_{\mathfrak{m}}(k)$. But, there is not much for us to prove here. Like in \cite{CDO},
$$| S_{\mathfrak{m}}(k)| = \frac{2^{r_u(k) + 1 + r_2(\mathrm{Cl}_{\mathfrak{m}}(k))}}{|\mathcal{O}_{k,\mathfrak{m}}/\mathcal{O}_{k,\mathfrak{m}}^2|},$$
where $r_u(k)$ is the rank of the unit group of $k$ and $r_2(\mathrm{Cl}_\mathfrak{m}(k))$ is the 2-rank of the ray class group of $k$ with modulus $\mathfrak{m}.$

Before proving Theorem \ref{local_CDO}, we will discuss the differences between our proof and the proof of Theorem 1.1 in \cite{CDO}. Lemma 3.3 from \cite{CDO} shows that quadratic extensions of any number field $k$ are in bijection with pairs $(\mathfrak{a},\bar{u})$ where $\mathfrak{a}$ is an integral, squarefree ideal with certain other conditions and $\bar{u}$ is a class in the $S(k).$ To construct $\Phi_{k,2}(C_2,s)$, they sum over integral, squarefree ideals $\mathfrak{a}$ with the appropriate conditions then sum over classes of $S(k)$ to pick up all pairs $(\mathfrak{a},\bar{u}).$ To handle ramification at 2, they look at all $\mathfrak{c}\mid 2$ and use the subgroup $S_{\mathfrak{c}^2}(k)$ to determine how many classes in $S(k)$ correspond to extensions with relative discriminant $\mathfrak{d}_{K,k} = 4\mathfrak{a}/\mathfrak{c}^2.$

In our proof we will construct the Dirichlet series by using the pairs $(\mathfrak{a},\bar{u}).$ However, because we are concerned with every class of $k_{\mathfrak{p}}^{\times}/k_{\mathfrak{p}}^{\times 2}$ for every $\mathfrak{p} \mid 2,$ we will handle ramification differently. As such, we will use a different modulus $\mathfrak{m}$ to set up our subgroup $S_{\mathfrak{m}}(k).$ Because every element of $k_{\mathfrak{p}}$ can be written as $\pi^i u$ where $\pi$ is the uniformizer, $i \in \Z$, and $u \in \mathcal{O}_{k_\mathfrak{p}}$, then we see that $k_{\mathfrak{p}}^{\times}/k_{\mathfrak{p}}^{\times 2} \cong \{\pm 1\}\times\mathcal{O}_{k_\mathfrak{p}}^{\times}/\mathcal{O}_{k_\mathfrak{p}}^{\times 2}.$ Moreover, for any $\mathfrak{p} \mid 2$, $\mathcal{O}_{k_\mathfrak{p}}^{\times}/\mathcal{O}_{k_\mathfrak{p}}^{\times 2} \cong (\mathcal{O}_k/\mathfrak{p}^{2e(\mathfrak{p}\mid 2)+1})^{\times}/(\mathcal{O}_k/\mathfrak{p}^{2e(\mathfrak{p}\mid 2)+1})^{\times 2}$. Thus the modulus we want is $\mathfrak{m} = \prod_{\mathfrak{p}\mid 2} \mathfrak{p}^{2e(\mathfrak{p}|2) + 1}$ and the complete set of local conditions can be represented globally by $\{(\pi_{\mathfrak{c}}, \bar{\beta}) \,\,|\,\, \pi_{\mathfrak{c}} = \prod_{\mathfrak{p}\mid \mathfrak{c}} \pi_\mathfrak{p}, \bar{\beta} \in \mathcal{O}_{k,\mathfrak{m}}/\mathcal{O}_{k,\mathfrak{m}}^2\},$ where $\mathfrak{c}$ is an ideal dividing $\mathfrak{n}$ as in the statement of Theorem \ref{local_CDO}. With this in mind, we are ready to prove the theorem.

\begin{proof}[(Proof of Theorem \ref{local_CDO})]
Cohen, Diay-y-Diaz, and Olivier start their proof of Theorem 1.1 with a sum over all squarefree $\mathfrak{a} \subset \mathcal{O}_k$ such that there exists $\mathfrak{q} \subset \mathcal{O}_k$ where $\mathfrak{a}\mathfrak{q}^2$ is principal. Because we will be manually handling cases where $\pi_{\mathfrak{c}} \ne 1$ we will instead start with a sum over $\mathfrak{c}\mid \mathfrak{n}$, where $\mathfrak{n} = \prod_{p \mid 2}\mathfrak{p}$ and require that there exist a $\mathfrak{q} \subset \mathcal{O}_k$ such that $\mathfrak{caq}^2$ is principal and $(\mathfrak{a},2) = 1$ So, we start with
$$\Phi_{k,2}(C_2, s) = -1 + \sum_{\mathfrak{c}\mid \mathfrak{n}}\sum_{\substack{\mathfrak{a} \sqfree \\ \exists \mathfrak{q}, \mathfrak{caq}^2 = \alpha_0\mathcal{O}_k \\ (\mathfrak{a},2)=1}}\mathrm{N}\mathfrak{a}^{-s}S(\alpha_0,\mathfrak{a},\mathfrak{c})$$
with
$$S(\alpha_0,\mathfrak{a},\mathfrak{c}) = \sum_{\bar{u} \in S(k)}D_2(\alpha_0u)^{-s}.$$
Since we want to count by the classes of $k_{\mathfrak{p}}^\times /k_{\mathfrak{p}}^{\times2}$, we use the correspondence between $\alpha_0u$ and $(\pi_{\mathfrak{c}},\bar{\beta})$ and rewrite the above sum as
$$S(\alpha_0,\mathfrak{a},\mathfrak{c}) = \sum_{\bar{\beta} \in \mathcal{O}_{k,\mathfrak{m}}/\mathcal{O}_{k,\mathfrak{m}}^2} D_2(\pi_{\mathfrak{c}}\beta)^{-s}f(\bar{\beta},\mathfrak{a},\mathfrak{c})$$
where $\beta$ is some lift of $\bar{\beta}$ to $\mathcal{O}_{k,\mathfrak{m}}$ and
$$f(\bar{\beta},\mathfrak{a},\mathfrak{c}) = \sum_{\substack{\bar{u} \in S(k) \\ \overline{\pi_{\mathfrak{c}}^{-1}\alpha_0u} = \bar{\beta}}}1.$$

We use $\overline{\pi_{\mathfrak{c}}^{-1}\alpha_0u} = \bar{\beta}$ to mean that the class of $\pi_\mathfrak{c}\alpha_0u$ in $\mathcal{O}_{k,\mathfrak{m}}/\mathcal{O}_{k,\mathfrak{m}}^2$ is $\bar{\beta}.$ This is equivalent to there existing a solution to $x^2 \equiv \pi_{\mathfrak{c}}^{-1}\beta^{-1}\alpha_0u \mods{\mathfrak{m}}.$ With Lemma \ref{ray_class_square}, this is equivalent to $\pi_{\mathfrak{c}}^{-1}\beta^{-1}\mathfrak{c}\mathfrak{a}$ being a square in $\mathrm{Cl}_{\mathfrak{m}}(k).$ Now, Lemma 3.10 from \cite{CDO} says that $f(\bar{\beta},\mathfrak{a},\mathfrak{c}) = 0$ if $\pi_{\mathfrak{c}}^{-1}\beta^{-1}\mathfrak{c}\mathfrak{a}$ is not a square in $\mathrm{Cl}_{\mathfrak{m}}(k)$ and $|S_{\mathfrak{m}}(k)|$ if it is. Theorem 2.36 from \cite{Pete} implies that $|\mathcal{O}_{k,\mathfrak{m}}/\mathcal{O}_{k,\mathfrak{m}}^2| = \prod_{\mathfrak{p}\mid 2}2 \cdot \mathrm{N}\mathfrak{p}^{2(e|\mathfrak{p})}.$ Adding this in to the cardinality of $S_{\mathfrak{m}}(k),$ we get
\begin{align*}
    |S_{\mathfrak{m}}(k)| &= \frac{2^{r_u(k) + 1 + r_2(\mathrm{Cl}_{\mathfrak{m}}(k))}}{\prod_{\mathfrak{p}\mid 2}(2 \cdot \mathrm{N}\mathfrak{p}^{e(\mathfrak{p\mid 2})})} \\
    &= \frac{2^{r_u(k) + 1 + r_2(\mathrm{Cl}_{\mathfrak{m}}(k))}}{2^{g+[k:\Q]}} \\
    &= \frac{2^{r_2(\mathrm{Cl}_{\mathfrak{m}}(k))}}{2^{g+i(k)}}.
\end{align*}

Like in \cite{CDO} we will use orthogonality and sum over the quadratic characters of $\mathrm{Cl}_{\mathfrak{m}}(k)$ to pick out the ideals that are squares of classes. So we have

\begin{align*}
    \Phi_{k,2}(C_2, s) &= -1 + \sum_{\mathfrak{c}\mid \mathfrak{n}}\sum_{\substack{\mathfrak{a} \sqfree \\ \exists \mathfrak{q}, \mathfrak{caq}^2 = \alpha_0\mathcal{O}_k \\ (\mathfrak{a},2)=1}} \frac{S(\alpha_0,\mathfrak{a},\mathfrak{c})}{\mathrm{N}\mathfrak{a}^s} \\
    &= -1 + \sum_{\mathfrak{c}\mid \mathfrak{n}}\sum_{\bar{\beta} \in \mathcal{O}_{k,\mathfrak{m}}/\mathcal{O}_{k,\mathfrak{m}}^2}\frac{|S_{\mathfrak{m}}(k)|}{D_2(\pi_{\mathfrak{c}}\bar{\beta})^{s}} \frac{1}{2^{r_2(\mathrm{Cl}_{\mathfrak{m}}(k))}}\sum_{\chi} \sum_{\substack{\mathfrak{a} \sqfree \\ (\mathfrak{a},2)=1}} \frac{\chi(\pi_{\mathfrak{c}}^{-1}\beta^{-1}\mathfrak{ca})}{\mathrm{N}\mathfrak{a}^s} \\
    &= -1 + \frac{1}{2^{g + i(k)}}\sum_{\mathfrak{c}\mid \mathfrak{n}}\sum_{\bar{\beta} \in \mathcal{O}_{k,\mathfrak{m}}/\mathcal{O}_{k,\mathfrak{m}}^2}\frac{1}{D_2(\pi_{\mathfrak{c}}\bar{\beta})^{s}} \sum_{\chi} \chi(\pi_{\mathfrak{c}}^{-1}\beta^{-1}\mathfrak{c}) \sum_{\substack{\mathfrak{a} \sqfree \\ (\mathfrak{a},2)=1}} \frac{\chi(\mathfrak{a})}{\mathrm{N}\mathfrak{a}^s}.
\end{align*}
At this point, the next few steps are identical to the end of the proof for Theorem 1.1 in \cite{CDO} with
$$\sum_{\substack{\mathfrak{a} \sqfree \\ (\mathfrak{a},2)=1}} \frac{\chi(\mathfrak{a})}{\mathrm{N}\mathfrak{a}^s} = \frac{L_k(s,\chi)}{\zeta_k(2s)\prod_{\mathfrak{p}\mid 2}(1 - \mathrm{N}\mathfrak{p}^{-2s})},$$ thus proving the theorem.
\end{proof}

If we let $N_{k,C_2}(X)$ stand for the number of quadratic extensions $K/k$ with $D_{K/k} \le X$, then one may use contour integration on our Dirichlet series from Theorem \ref{local_CDO} to get the same result as their Corollary 1.2. To limit the count to only quadratic fields $K/k$ with a specific local condition at 2, simply limit the sum to the set of $(\pi_{\mathfrak{c}}\,\bar{\beta})$ in question. But, we can be more general than this and consider local specifications for other places of $k$.

Recalling our definition for an acceptable collection of local specifications $\Sigma_k$, let $S(\Sigma_k) = \{v \mid 2\} \cup \{v $ s.t. $\Sigma_{k,v}$ doesn't contain every degree 2 \'etale algebra of $k_v\}.$ We assume that $\Sigma_{k,v}$ contains every degree 2 \'etale algebra when $v$ is a complex infinite place as there is only one. This set contains all of the places with local information we need to write out $\Phi_{k,2}(\Sigma_k; C_2, s)$ with full detail. Each $K_v \in \Sigma_{k,v}$ corresponds to a class of $(\mathcal{O}_{k_{v}}/(\pi_v))/(\mathcal{O}_{k_{v}}/(\pi_v))^2$ for finite places and $\sigma_v(x) > 0$ or $< 0$ where $K_v = k_v(\sqrt{x})$ for the real infinite places. Using the weak approximation theorem, this extends to a set $B$ containing elements of the form $(\pi,\bar{\beta})$ where $\pi$ is a product (potentially empty) of uniformizers $\pi_v \in k_v$ for finite $v \in S(\Sigma_k)$ and $\bar{\beta} \in \mathcal{O}_{k,\mathfrak{m}}/\mathcal{O}_{k,\mathfrak{m}}^2$ for a modulus $\mathfrak{m}$ with $\mathfrak{m}_0 = \prod_{\mathfrak{p}\mid 2}\mathfrak{p}^{2e+1}\prod_{\substack{v \text{ odd}\\ v \in S(\Sigma_k)}} \mathfrak{p}_v$ and $\mathfrak{m}_{\infty} = \prod_{\substack{v \text{ real}\\ v \in S(\Sigma_k)}}\sigma_v.$ We require that each $\pi_\mathfrak{p}$ be $1 \mods{q}$ for every other $\mathfrak{q} \mid \mathfrak{m}_0$ and $\sigma_v(\pi_\mathfrak{p}) > 0$ for every $\sigma_v\mid \mathfrak{m}_\infty.$ We will write these elements as tuples $(\pi_{\mathfrak{c}}, \bar{\beta}) \in B$ where $\mathfrak{c}$ is some squarefree ideal divisor of $\mathfrak{m}_0$ and $\pi_{\mathfrak{c}}$ is defined as it was earlier.

\begin{theorem}\label{local_dirichlet_series}
Let $k$ be a number field, $\Sigma_k$ be an acceptable collection of local specifications for $k$, and $B$ be the corresponding set of global conditions for $K \in \mathcal{K}(\Sigma_k).$ Let $\mathfrak{m}=\mathfrak{m}_0\mathfrak{m}_\infty$ be the corresponding modulus.  Then,
$$\Phi_{k,2}(\Sigma_k; C_2, s) = -\delta_{(1,\bar{1}) \in B} + \frac{1}{2^{|S(\Sigma_k)| + i(k)}\zeta_k(2s)} \prod_{\mathfrak{p} \mid \mathfrak{m}_0}\left(1 - \mathrm{N}\mathfrak{p}^{-2s}\right)^{-1} \sum_{(\pi_{\mathfrak{c}},\bar{\beta}) \in B}D_{\mathfrak{m}_0}(\pi_\mathfrak{c} \beta)^{-s} \sum_{\chi}\chi(\pi_\mathfrak{c}^{-1}\beta^{-1} \mathfrak{c})L_k(s,
\chi),$$
where $D_{\mathfrak{m}_0}(x) = \prod_{\mathfrak{p}\mid \mathfrak{m}_0}D_{k_\mathfrak{p}(\sqrt{x})/k_\mathfrak{p}}$ and $\delta_{(1,\bar{1}) \in B}$ is $1$ if $(1,\bar{1}) \in B$ and $0$ otherwise.
\end{theorem}

\begin{proof}
The proof is almost identical to that of Theorem \ref{local_CDO} save for few details. When choosing the lift $\beta$ of $\bar{\beta}$, for each $\sigma_v \mid \mathfrak{m}_\infty$, we must require that $\sigma_v(\beta)$ be either greater than or less than $0$ depending on if the corresponding $K_v \in \Sigma_{k,v}$ is $\R^2$ or $\C,$ respectively. Also we now  need to evaluate the cardinality of $\mathcal{O}_{k,\mathfrak{m}}/\mathcal{O}_{k,\mathfrak{m}}^2$ for a general modulus $\mathfrak{m}.$ For this, we note that for any odd prime $\mathfrak{p}$, $\mathcal{O}_{k_{\mathfrak{p}}}^{\times}/\mathcal{O}_{k_{\mathfrak{p}}}^{\times2} \cong (\mathcal{O}_k/\mathfrak{p})^{\times}/(\mathcal{O}_k/\mathfrak{p})^{\times 2}$ and $|(\mathcal{O}_k/\mathfrak{p})^{\times}/(\mathcal{O}_k/\mathfrak{p})^{\times 2}| = 2.$ Also, for any real infinite place $v$, $|\mathcal{O}_{k_{v}}^{\times}/\mathcal{O}_{k_{v}}^{\times2}| =2.$
\end{proof}

We now have everything we need to prove Theorem \ref{local_asymptotics}.
\begin{proof}
To prove the theorem, we will apply Perron's formula on $\Phi_{k,2}(\Sigma_k; C_2,s)$ as given by Theorem \ref{local_dirichlet_series} as this will give us $N_{k,C_2}(\Sigma_k; X).$ Let $\Sigma_k$ be an acceptable collection of local specifications for $k$. We take the modulus $\mathfrak{m}$ to be defined as in the statement of Theorem \ref{local_asymptotics}. 

We can exclude $\mathfrak{r}_1$ (resp. $\mathfrak{u}_1$) from the modulus $\mathfrak{m}$ and the set $B$ because we don't need to pick subset of the classes in $\mathcal{O}_{k,\mathfrak{r}_1}^{\times}\mathcal{O}_{k,\mathfrak{r}_1}^{\times2}$ (resp. $\mathcal{O}_{k,\mathfrak{u}_1}^{\times}\mathcal{O}_{k,\mathfrak{u}_1}^{\times2}$). To count these extensions we modify the sum over integral, squarefree ideals to be

$$\frac{1}{\mathrm{N}\mathfrak{r}_1^2} \sum_{\substack{\mathfrak{a}\sqfree \\ (\mathfrak{a},\mathfrak{m}_0\mathfrak{r}_1)=1}} \frac{\chi(\mathfrak{a})}{\mathrm{N}\mathfrak{a}^2}.$$
This way, we are using a character $\chi$ from $\mathrm{Cl}_\mathfrak{m}(k)$ but we can modify the Euler products for $L_k(s,\chi)$ and $\zeta_k(2s)$ such that they exclude primes dividing $\mathrm{r}_1.$ A similar change to the sum can be made for the primes dividing $\mathfrak{u}_1.$

Now, note that all of the $L$-functions in the sum over quadratic characters of the ray class group are holomorphic on the entire complex plane with the exception of the trivial character, which has a pole at $s=1$. So, $\Phi_{k, 2}(\Sigma_k; C_2, s)$ only has a pole at $s=1$ for $\mathrm{Re}(s)>1/2.$

Using the standard convexity bounds on Hecke characters (e.g. \cite[p.~142]{IK}), we get that
$$N_{k,2}(C_2, \Sigma_k, X) = \res\Phi_{k,2}(\Sigma_k; C_2, s)X + O_{n,\epsilon}\left(\#\mathrm{Cl}_{\mathfrak{m}}(k)[2]\frac{X^{\frac{n+2}{n+4} + \epsilon}|D_K|^{\frac{1}{n+4}+\epsilon}\mathrm{N}\mathfrak{u}_1^{\epsilon}\mathrm{N}\mathfrak{u}_2^{\frac{1}{n+4} + \epsilon}}{\mathfrak{r}_1^{\frac{n+2}{n+4} - \epsilon}\mathfrak{r}_1^{\frac{n+1}{n+4} - \epsilon}}\right).$$
Here, our modifications to the sum for $\mathfrak{r}_1$ and $\mathfrak{u}_1$ give a better error estimate because they are unaffected by the modulus of $\chi$ in the convexity bound.

Examining the residue of the Dirichlet function, we see that
\begin{align*}
    \res \Phi_{k,2}(\Sigma_k;C_2,s) &= \frac{\res L_k(s,\chi_0)}{2^{|S(\Sigma_k)|+i(k)}\zeta_k(2)} \prod_{\mathfrak{p}\mid \mathfrak{r}_1}\left(1 + \mathrm{N}\mathfrak{p}\right)^{-1}\prod_{\mathfrak{p}\mid \mathfrak{u}_1}\left(1 + \frac{1}{\mathrm{N}\mathfrak{p}}\right)^{-1}\prod_{\mathfrak{p}\mid \mathfrak{m}_0}\left(1 - \mathrm{N}\mathfrak{p}^{-2}\right)^{-1}\sum_{(\pi_{\mathfrak{c}},\bar{\beta}) \in B}\frac{1}{D_{\mathfrak{m}_0}(\pi_\mathfrak{c}\beta)} \\
    &= \frac{\res \zeta_k(s)}{2^{|S(\Sigma_k)|+i(k)}\zeta_k(2)}\prod_{\mathfrak{p}\mid \mathfrak{r}_1}\left(1 + \mathrm{N}\mathfrak{p}\right)^{-1}\prod_{\mathfrak{p} \mid \mathfrak{m}_0\mathfrak{u}_1} \left(1 + \frac{1}{\mathrm{N}\mathfrak{p}}\right)^{-1}\sum_{(\pi_{\mathfrak{c}},\bar{\beta}) \in B}\frac{1}{D_{\mathfrak{m}_0}(\pi_\mathfrak{c}\bar{\beta})} \\
    & = \frac{\res \zeta_k(s)}{2^{|S(\Sigma_k)|+i(k)}\zeta_k(2)}\prod_{\sigma_v \mid \infty}\left(\sum_{K_v \in \Sigma_{k,v}} 1\right)\prod_{\mathfrak{p} \mid \mathfrak{m}_0} \left(\left(1 + \frac{1}{\mathrm{N}\mathfrak{p}}\right)^{-1}\sum_{K_{\mathfrak{p}} \in \Sigma_{k,\mathfrak{p}}}\frac{1}{D_{K_{\mathfrak{p}}/k_{\mathfrak{p}}}} \right) \\
    & = \frac{\res \zeta_k(s)}{\zeta_k(2)}\prod_{\sigma_v \mid \infty}\left(\sum_{K_{v} \in \Sigma_{k,v}}\frac{1}{\Aut{K_{v}/k_{v}}}\right)\prod_{\mathfrak{p}}\left(\left(1 + \frac{1}{\mathrm{N}\mathfrak{p}}\right)^{-1}\sum_{K_{\mathfrak{p}} \in \Sigma_{k,\mathfrak{p}}} \frac{1}{\Aut{K_{\mathfrak{p}}/k_{\mathfrak{p}}}D_{K_{\mathfrak{p}}/k_{\mathfrak{p}}}} \right).
\end{align*}
We are able to extend the product to all finite primes by noting that
$$\left(1 + \frac{1}{\mathrm{N}\mathfrak{p}}\right)^{-1}\sum_{K_{\mathfrak{p}} \in \Sigma_{k,\mathfrak{p}}} \frac{1}{\Aut{K_{\mathfrak{p}}/k_{\mathfrak{p}}}D_{K_{\mathfrak{p}}/k_{\mathfrak{p}}}} = 1$$
whenever $\Sigma_{k,\mathfrak{p}}$ contains every degree 2 \'etale algebra of $k_{\mathfrak{p}}.$
\end{proof}

We now have all the tools we need and can move on to the proof of the main theorem.

\section{Proving Theorems \ref{main_thm} and \ref{local_specs}}
We will first prove the main theorem and then present an outline for proving Theorem \ref{local_specs}, as it is very similar. We will prove Theorem \ref{main_thm} by considering the different pieces of (\ref{the_sum}) separately beginning with the double sum
$$\sum_{\substack{[K:\Q]=2\\ |D_K| < X^{1/2}}}\sum_{\substack{[L:K]=2 \\ D_{L/K} < X/|D_K|}} 1.$$
We cannot simply use Theorem 5.3 from \cite{ASVW} as the error term is too large for our purposes. The key for us will be to prove our own version of Theorem 2 from \cite{ASVW}.

\begin{lemma}\label{theorem_2}
For any $X \ge 1$,
$$\sum_{\substack{[K:\Q]=2 \\ 0 < D_K < X}} \frac{1}{|D_K|}\cdot \frac{L(1,K/\Q)}{L(2,K/\Q)} = (\log X + 1)\frac{\zeta(2)}{2} \prod_{p}\left(1 - \frac{1}{p^2} - \frac{2}{p^3} + \frac{2}{p^4}\right) + c^+ + O_{\epsilon}(X^{-5/18 + \epsilon});$$
$$\sum_{\substack{[K:\Q]=2 \\  -X< D_K < 0}} \frac{1}{|D_K|}\cdot \frac{L(1,K/\Q)}{L(2,K/\Q)} = (\log X + 1)\frac{\zeta(2)}{2} \prod_{p}\left(1 - \frac{1}{p^2} - \frac{2}{p^3} + \frac{2}{p^4}\right) + c^- + O_{\epsilon}(X^{-5/18 + \epsilon}),$$
for some constants $c^+$ and $c^-$.
\end{lemma}

The proof of this lemma follows a similar outline to Theorem 2 of \cite{ASVW}. Where they approximate $L(1, K/\Q)$ using a finite sum, we will use a smooth approximation instead.
\begin{lemma}\label{smooth_approx}
Let $L(s,K/\Q)$ be the Dirichlet $L$-function attached to the quadratic extension $K/\Q.$ For any $N>1$,
$$L(1, K/\Q) = \sum_{n \ge 1}\frac{\chi_K(n)e^{-n/N}}{n} + O_{\epsilon}\left(\frac{D_K^{1/6 + \epsilon}}{N^{1/2}}\right),$$
where $\chi_K$ is the Kronecker character for $K/\Q.$
\end{lemma}
\begin{proof}
Using an inverse Mellin transform similar to Perron's formula, we have that for any $c > 0$
\begin{align*}
    \sum_{n \ge 1}\frac{\chi_K(n)e^{-n/N}}{n}& = \frac{1}{2 \pi i}\int_{(c)} L(s+1, K/\Q)N^{s}\Gamma(s)\,ds \\
    & = L(1, K/\Q) + \frac{1}{2\pi i}\int_{(-1/2)} L(s + 1, K/\Q)N^{s}\Gamma(s)\,ds,
\end{align*}
where $\int_{(c)}$ indicates the contour integral from $c-i\infty$ to $c+i\infty$ and $\int_{(-1/2)}$ indicates a contour integral where part of the contour is pushed to $\Re(s)=-1/2.$
The lemma follows from the subconvexity bound $L(1/2 + it, K/\Q) \ll_{\epsilon} (D_K(1+t))^{1/6 + \epsilon}$ of Petrow and Young \cite{PY}.
\end{proof}

Next, we look at sums over $\chi_K(n)$ as $K$ varies and see that the largest contributions come from $n$ being square.

\begin{lemma}\label{sum_kronecker}
For any integer $n \ge 1$ and number $X \ge 1$, if $n$ is not a square, 
$$\sum_{\substack{[K:\Q] = 2 \\ X < D_K < 2X}} \chi_K(n) \ll X^{1/2}n^{1/4 + \epsilon},$$
and
$$\sum_{\substack{[K:\Q] = 2 \\ X < D_K < 2X}} \chi_K(n) = \frac{X}{2\zeta(2)} \prod_{p \mid n} \frac{p}{p + 1} + O(X^{1/2}n^{\epsilon})$$
if $n$ is a square.
\end{lemma}
\begin{proof}
When $n$ is not a square, we rewrite $\chi_K(n)$ as $\left(\frac{D_K}{n}\right)$ use the inclusion-exclusion principle so we can sum over squarefree integers. Thus our sum is
$$\sum_{a < \sqrt{2X}}\mu(a) \sum_{\frac{X}{a^2} < d < \frac{2X}{a^2}}\left(\frac{a^2d}{n}\right).$$
We now split the sum over $a$ into two parts using an auxillarly parameter $T$. We will treat the first part with Polya-Vinogradov and the second with a trivial estimate.

\begin{align*}
    \sum_{a < \sqrt{2X}}\mu(a) \sum_{\frac{X}{a^2} < d < \frac{2X}{a^2}}\left(\frac{a^2d}{n}\right) &= \sum_{a < T}\mu(a) \sum_{\frac{X}{a^2} < d < \frac{2X}{a^2}}\left(\frac{a^2d}{n}\right) + \sum_{T \le a < \sqrt{2X}}\mu(a)\sum_{\frac{X}{a^2} < d < \frac{2X}{a^2}}\left(\frac{a^2d}{n}\right) \\
    &= \sum_{a < T} O(n^{1/2 + \epsilon}) + \sum_{T \le a < \sqrt{2X}} O(X/a^2) \\
    &= O(Tn^{1/2+\epsilon} + X/T).
\end{align*}
Choosing $T = X^{1/2}/n^{1/4}$ yields the first part of the lemma.

When $n$ is a square, the result follows from Theorem \ref{local_asymptotics} but using results from counting squarefree integers to obtain the error term. See \cite[Lemma~2.17]{MV} for a reference.
\end{proof}

We then combine these lemmas to analyze the sum of $\frac{L(1,K/\Q)}{L(2,K/\Q)}$ for $D_K$ in the range $X$ to $2X$.

\begin{lemma}\label{smooth_l}
For any $X \ge 1$,
$$\sum_{\substack{[K:\Q]=2 \\X < D_K < 2X}} \frac{L(1,K/\Q)}{L(2,K/\Q)} = \frac{1}{2\zeta(2)}X \prod_p \left(1 + \frac{1}{(p+1)^2}\right) + O_{\epsilon}(X^{13/18 + \epsilon}).$$
\end{lemma}
\begin{proof}
We start be rewriting $1/L(2,K/\Q)$ as the Dirichlet series
$$\sum_{m \ge 1}\frac{\mu(m)\chi_K(m)}{m^2}$$ and combining this with Lemma \ref{smooth_approx}. For $N > X^{1/3}$, we see that
$$\frac{L(1, K/\Q)}{L(2,K/\Q)} = \sum_{m \ge 1}\sum_{n \ge 1}\frac{\mu(m)\chi_K(mn)e^{-n/N}}{m^2n} + O_{\epsilon}\left(\frac{D_K^{1/6+\epsilon}}{N^{1/2}}\right).$$
Therefore, summing $D_K$ in the range $X$ to $2X$ gives
$$\sum_{\substack{[K:\Q]=2 \\ X < D_K < 2X}}\frac{L(1, K/\Q)}{L(2,K/\Q)} =  \sum_{m \ge 1}\sum_{n \ge 1}\frac{\mu(m)e^{-n/N}}{m^2n}\sum_{\substack{[K:\Q]=2 \\ X < D_K < 2X}}\chi_K(mn) + O_{\epsilon}\left(\frac{X^{7/6 + \epsilon}}{N^{1/2}}\right).$$
By Lemma \ref{sum_kronecker}, the inner sum is negligible unless $mn$ is a square. Because $m$ must also be squarefree, this can only happen when $m$ is the squarefree part of $n$. If we let $\psi(n)$ be the multiplicative function such that $\psi(p^k) = p/(p +1)$, we can eliminate the sum over $m$ because
$$\sum_{m \ge 1}\sum_{n \ge 1}\frac{\mu(m)e^{-n/N}}{m^2n}\sum_{\substack{[K:\Q]=2 \\ X < D_K < 2X}}\chi_K(mn) = \frac{X}{2\zeta(2)}\sum_{n \ge 1}\frac{\mu(m)\psi(n)e^{-n/N}}{m^2n} + O(X^{1/2+\epsilon}N^{1/4 + \epsilon}).$$
We can simplify this even further.
\begin{align*}
    \frac{X}{2\zeta(2)}\sum_{n \ge 1}\frac{\mu(m)\psi(n)e^{-n/N}}{m^2n} + O(X^{1/2+\epsilon}N^{1/4 + \epsilon}) &= \frac{X}{2\zeta(2)}\sum_{n \ge 1}\frac{\mu(m)\psi(n)}{m^2n} + O(X/N + X^{1/2+\epsilon}N^{1/4 + \epsilon}) \\
    &= \frac{X}{2\zeta(2)}\prod_p \left(1 + \frac{1}{(p + 1)^2}\right) + O(X/N + X^{1/2+\epsilon}N^{1/4 + \epsilon}).
\end{align*}
Though the penultimate sum does not initially look like it should yield an Euler product of this form, analyzing it shows that there are geometric sub-series that, when combined, lead to the final product. Now, to optimize both error terms, we choose $N = X^{8/9}.$
\end{proof}
We can now prove Lemma \ref{theorem_2} by using the above lemmas and partial summation.
\begin{proof}[(Proof of Lemma \ref{theorem_2})]
Let $S^+(X)$ denote $\sum_{\substack{[K:\Q]=2 \\0 < D_K < X}} \frac{L(1,K/\Q)}{L(2,K/\Q)}.$ Then
\begin{align*}
    \sum_{\substack{[K:\Q]=2 \\ 0 < D_K < X}} \frac{1}{|D_K|}\cdot \frac{L(1,K/\Q)}{L(2,K/\Q)} &= \frac{S^+(t)}{t}\Big|_{1^{-}}^{X} + \int_{1^{-}}^{X} \frac{S^+(t)}{t^2}dt \\
    &=\frac{1}{2\zeta(2)}\prod_{p}\left(1 + \frac{1}{(p+1)^2}\right)(\log X + 1) + \int_{1^{-}}^{X} \frac{E^+(t)}{t^2}dt +  O_{\epsilon}(X^{-5/18 + \epsilon}) \\
    &= \frac{1}{2\zeta(2)}\log(X)\prod_{p}\left(1 + \frac{1}{(p+1)^2}\right) + \int_{1^{-}}^{\infty} \frac{E^+(t)}{t^2}dt + O_{\epsilon}(X^{-5/18 + \epsilon}).
\end{align*}
Above, $E^+(X) = S^+(X) - \frac{1}{2\zeta(2)}X\prod_{p}\left(1 + \frac{1}{(p+1)^2}\right)$ with $\int_{1^-}^X \frac{E^+(t)}{t^2}dt = \int_{1^{-}}^{\infty} \frac{E^+(t)}{t^2}dt - \int_{X}^{\infty} \frac{E^+(t)}{t^2}dt.$ This is justified because $E(X) = O(X^{13/18 + \epsilon})$, so the integral $\int_{1^{-}}^{\infty} \frac{E^+(t)}{t^2}dt$ converges. With this, we define
\begin{equation}
    c^+ = \int_{1^{-}}^\infty \frac{E^+(t)}{t^2}dt.
\end{equation}
Defining $S^-(X)$ and $E^-(X)$ similarly, we also define
\begin{equation}
    c^- = \int_{1^{-}}^\infty \frac{E^-(t)}{t^2}dt.
\end{equation}

At this point, note that $$\frac{1}{\zeta(2)}\prod_{p}\left(1 + \frac{1}{(p+1)^2}\right) = \zeta(2)\prod_p\left(1 -\frac{1}{p^2} - \frac{2}{p^3} + \frac{2}{p^4}\right).$$
\end{proof}

We can now use Theorem \ref{local_asymptotics} on the inner sum and Lemma \ref{theorem_2} on the outer sum to get
\begin{equation}\label{first_sum}
    \sum_{\substack{[K:\Q]=2\\ |D_K| < X^{1/2}}}\sum_{\substack{[L:K]=2 \\ D_{L/K} < X/|D_K|}} 1 = X \left(\frac{1}{2}\log X + 1\right) \frac{3}{4} \prod_p \left(1 - \frac{1}{p^2} - \frac{2}{p^3} + \frac{2}{p^4} \right) + cX + O_{\epsilon}(X^{11/12 + \epsilon}),
\end{equation}
where
\begin{equation}\label{c_defn}
c = c^+ + \frac{1}{2}c^-.
\end{equation}
 
\subsection{The Second Sum}
We will now turn our attention to the second sum from (\ref{the_sum}), namely
$$\sum_{q < X^{1/4}}\sum_{\substack{[K:\Q]=2\\ |D_K| < X^{1/2}/q^2}}\sum_{\substack{[L:K]=2 \\ D_{L/K} < X^{1/2}q^2 \\ J(L) = q^2}} 1.$$

Now that we can count quadratic extensions $L/K$ with $J_2(L) = 2^{2i}$ for some $i = 0, 1, 2,$ or $3$, we will expand $q$ to get $J(L)= 2^{2i}d^2$. So, our second sum is now
\begin{equation}\label{four_sum}
    \sum_{\substack{q=2^i d < X^{1/4}\\ d \text{ odd,}\sqfree \\i \in \{0, 1, 2, 3\}}}\sum_{\substack{[K:\Q]=2\\ |D_K| < X^{1/2}/q^2}}\sum_{\substack{[L:K]=2 \\ D_{L/K} < X^{1/2}q^2 \\ J(L) = q^2}} 1.
\end{equation}

For the odd part of $J(L)$, we will require that $(D_K, d) = 1$, and for $J_2(L)$, we will require that $D_K$ fall into specific congruence classes mod $8$ corresponding to Tables 1-4. To limits ourselves to extensions $L/K$ such that the odd part of $J(L)$ is exactly $d^2,$ we will use the inclusion/exclusion principle. Thus, the innermost sum of (\ref{four_sum}) becomes
\begin{equation}\label{before_e}
    \sum_{\substack{[L:K]=2 \\ D_{L/K} < X^{1/2}q^2 \\ J(L) = q^2}} 1 = \sum_{\substack{e < X^{1/4}\\(e, |D_K|2d)=1}} \mu(e) \sum_{\substack{[L:K]=2\\D_{L/K} < X^{1/2}q^2\\d^2e^2 \mid D_{L/K}\\J_2(L)=2^{2i}}} 1.
\end{equation}
Before we apply Theorem \ref{local_asymptotics}, we construct an acceptable collection of local conditions $\Sigma$ for which $\Sigma_2$ only contains pairs $(L_2,K_2)$ where $J_2(L) = 2^{2i}$ and $\Sigma_p$ contains pairs $(L_p,K_p)$ for all $p \mid de$  such that $J_p(L)=p^2.$ $\Sigma$ can be restricted to an acceptable collection $\Sigma_K$ of local conditions for $K$. As we vary $K$ later on in the argument, we will pick up every pair $(L_p,K_p) \in \Sigma_p$ for each $p \mid 2de.$ Now, we apply the theorem on $\Sigma_K$ and get

$$X^{1/2}\frac{q^2L(1,K/\Q)}{2^{i(K)}\zeta(2)L(2,K/\Q)}\prod_{\mathfrak{p}\mid 2de}\left(\left(1 + \frac{1}{\mathrm{N}\mathfrak{p}}\right)^{-1}\sum_{L_{\mathfrak{p}} \in \Sigma_{K,{\mathfrak{p}}}}\frac{1}{\Aut{L_{\mathfrak{p}}/K_{\mathfrak{p}}}D_{L_{\mathfrak{p}}/K_{\mathfrak{p}}}}\right) + O_{\epsilon}\left(\frac{X^{1/3+\epsilon}q^{\epsilon}|D_K|^{1/6+\epsilon}}{e^{4/3 + \epsilon}}\right).$$

We will pull everything but the product over primes dividing $e$ in front of our sum over $e$. For those remaining primes, because $\Sigma_{K,\mathfrak{p}}$ contains only the $L_{\mathfrak{p}}$ that ramify we see that 
\begin{equation} \label{simple_e}
\prod_{\mathfrak{p}\mid e}\left(1 + \frac{1}{\mathrm{N}\mathfrak{p}}\right)^{-1}\sum_{L_{\mathfrak{p}} \in \Sigma_{K,{\mathfrak{p}}}}\frac{1}{\Aut{L_{\mathfrak{p}}/K_{\mathfrak{p}}}D_{L_{\mathfrak{p}}/K_{\mathfrak{p}}}} = \prod_{p \mid e} \frac{1}{m(p)},
\end{equation}
where $m(p) = 1/(1 + p)^2$ if $p$ splits in $K$ and $1/(1 + p^2)$ if $p$ is inert in $K$. So, the sum over $e$ is 
\begin{align*}
    \sum_{\substack{e < X^{1/4}\\(e, |D_K|2d)=1}}\left( \mu(e)\prod_{p\mid e}\frac{1}{m(p)}  + O_{\epsilon}\left(\frac{X^{1/3+\epsilon}q^{\epsilon}|D_K|^{1/6+\epsilon}}{e^{4/3 + \epsilon}}\right)\right)  =  \prod_{p \nmid D_K 2d}\left(1 - \frac{1}{m(p)}\right) + O_{\epsilon}(X^{1/3 + \epsilon} q^{\epsilon}|D_K|^{1/6 + \epsilon}).
\end{align*}

Before we turn to the sum over $K,$ we will deal with some of the constants we have accrued. First, we let

$$\mu(\Sigma_{K,2^{2i}}) = \prod_{\mathfrak{p}\mid 2}\sum_{L_{\mathfrak{p}}\in \Sigma_{K,\mathfrak{p}}} \frac{1}{\Aut{L_{\mathfrak{p}}/K_{\mathfrak{p}}}D_{L_{\mathfrak{p}}/K_{\mathfrak{p}}}}.$$
Additionaly, we see that if pull out $p=2$ from the Euler product for $\frac{L(1,K/\Q)}{L(2,K/\Q)}$ we get
\begin{equation}\label{simple_p}
    \left(\frac{1 - \chi_K(p)/p^2}{1 - \chi_K(p)/p}\right)\prod_{\mathfrak{p} \mid 2}\left(1 + \frac{1}{\mathrm{N}\mathfrak{p}}\right)^{-1} = \left(1 + \frac{1}{p}\right)^{-1} 
\end{equation}
regardless of if $2$ splits, is inert, or ramifies in $K.$

Similarly, we also consider the factors of $\frac{L(1,K/\Q)}{L(2,K/\Q)}$ at odd primes $p$. For the primes dividing $d$, we can use (\ref{simple_e}) and get
\begin{equation}\label{dividing_q}
    \left(\frac{1}{m(p)}\right)\left(\frac{1 - \chi_K(p)/p^2}{1 - \chi_K(p)/p}\right) = \frac{1}{p(p+1)}.
\end{equation}
For the remaining primes not dividing $D_K$,
\begin{equation}\label{dividing_e}
\left(1 - \frac{1}{m(p)}\right)\left(\frac{1 - \chi_K(p)/p^2}{1 - \chi_K(p)/p}\right) = 1 + \frac{\chi_K(p)}{p+1}.
\end{equation}

Applying all of this, the sum over $K$ is
\begin{equation}\label{after_e}
    X^{1/2}\frac{2q^2}{3\zeta(2)}\prod_{p \mid d}\frac{1}{p(p+1)} \sum_{\substack{[K:\Q]=2\\|D_K|<X^{1/2}/q^2 \\(D_K,d)=1\\K_2 \in \Sigma_{\Q,2^{2i}}}} \frac{\mu(\Sigma_{K,2^{2i}})}{2^{i(K)}}\prod_{p \nmid 2d}\left(1 + \frac{\chi_K(p)}{p+1}\right) + O_{\epsilon}\left(\frac{X^{11/12 + \epsilon}}{q^{7/3 - \epsilon}}\right),
\end{equation}
where $\Sigma_{\Q, 2^{2i}}$ is an acceptable set of local conditions over $\Q$ that is complete at every odd prime and at $p=2$ contains only quadratic extensions $K$ that are quadratic subfields for quartic fields $L$ with $J_2(L) = 2^{2i}.$

To analyze the sum, we first rewrite the Euler product as the sum

\begin{equation}\label{to_swap}
    \sum_{\substack{[K:\Q]=2\\|D_K|<X^{1/2}/q^2\\(D_K,d)=1\\K_2 \in \Sigma_{\Q,2^{2i}}}} \frac{\mu(\Sigma_{K,2^{2i}})}{2^{i(K)}} \sum_{\substack{n \ge 1 \\ (n,2d)=1\\n \sqfree}}\frac{\chi_K(n)}{f(n)},
\end{equation}
where $f(n) = \prod_{p \mid n} (p +1).$ In a series of lemmas, we will consider how to estimate the inner sum of (\ref{to_swap}) and swap the order of summation to arrive at our conclusion. We again use the inverse Mellin transfrom from Lemma \ref{smooth_approx}. To do this, we consider the following function.
$$ F(s, \chi_K) = \sum_{\substack{n \ge 1\\(n,d)=1}} \frac{\mu^2(n)\chi_K(n)\psi(n)}{n^s},$$
where $\psi(n) = \prod_{p \mid n}\frac{p}{p+1}.$ Note that $F(1,\chi_K)$ is equivalent to the inner sum of (\ref{to_swap}).

\begin{lemma}\label{smooth_perron}
Let $F(s, \chi_K)$ be the function given above. Then for any $N > 1,$
$$ F(1, \chi_K) = \sum_{\substack{n \ge 1\\ (n,d)=1\\n \sqfree}} \frac{\chi_K(n)e^{-n/N}}{f(n)} + O_{\epsilon}\left(\frac{|D_K|^{1/6 + \epsilon}}{N^{1/2}}\right).$$
\end{lemma}

\begin{proof}
If we note that
\begin{align*}
    F(s, \chi_K) &= \prod_{p\nmid d}\left(1 + \frac{\chi_K(p)}{p^{s-1}(p+1)}\right) \\
    &= L(s,\chi_K)\prod_{p \mid d}\left(1 + \frac{\chi_K(p)}{p^s}\right)\prod_{p \nmid d}\left(1 - \frac{\chi_K(p)}{p^s(p+1)} - \frac{1}{p^{2s-1}(p+1)}\right),
\end{align*}
and the rightmost product is absolutely convergent for $\Re(s+1) > -1/2$, then the proof for this lemma is almost identical to Lemma \ref{smooth_approx}.
\end{proof}

In addition to this, we also need two more lemmas also very similar to Lemmas \ref{sum_kronecker} and \ref{smooth_l}. In the lemmas below, because $\mu(\Sigma_{K,2^{2i}})$ depends on the class that $K_2$ corresponds to in $\Q_2/\Q_2^{\times 2}$, we sum over $K$ such that $D_K \equiv a \mod{8}$ rather than $K_2 \in \Sigma_{\Q, 2^{2i}}$, where the value of $a$ corresponds to the congruence conditions from Tables 1-4. Then, when we use the lemmas to analyze (\ref{to_swap}), we will combine the separate congruence conditions to get the sum over $K_2 \in \Sigma_{\Q, 2^{2i}}.$

\begin{lemma}\label{sum_chi}
For any squarefree integers $n, d \ge 1$ with $(n,d)=1, X \ge 1,$ and $a \in \{0, 1, 4, 5\}$ we have
$$\sum_{\substack{[K:\Q]=2\\X < D_K < 2X\\(D_K,d)=1\\ D_K \equiv a \mod{8}}} \chi_K(n) \ll_{\epsilon} X^{1/2+\epsilon}(nd)^{1/4 + \epsilon}$$ if $n \ne 1$, and
$$\sum_{\substack{[K:\Q]=2\\X < D_K < 2X\\(D_K,d)=1\\ D_K \equiv a \mod{8}}} \chi_K(n) = X\frac{\mu(\Sigma_{\Q, 2^{2i},a})}{3\zeta(2)}\prod_{p \mid d}\frac{p}{p+1} + O_{\epsilon}(X^{1/2}d^{\epsilon})$$ if $n=1,$
where 
$$\mu(\Sigma_{\Q, 2^{2i},a}) = \sum_{\substack{K_2 \in \Sigma_{\Q, 2^{2i}} \\ D_K \equiv a \mod{8}}}\frac{1}{\Aut{K_2/\Q_2}D_{K_2/\Q_2}}.$$
\end{lemma}

\begin{proof}
The proof for this is nearly identical to Lemma \ref{sum_kronecker} except that we consider $\chi_K(n)$ as a character with conductor $nd$.
\end{proof}

Lastly, the clone of Lemma \ref{smooth_l}.

\begin{lemma}\label{diadic_sum}
For any $X \ge 1$, we have
$$\sum_{\substack{[K:\Q]=2\\X < D_K < 2X\\(D_K,d)=1\\D_K \equiv a \mod{8}}}\sum_{\substack{n\ge1 \\ (n,2d)=1 \\ n\, \sqfree}} \frac{\chi_K(n)}{f(n)} = X\frac{\mu(\Sigma_{\Q, 2^{2i}, a})}{3\zeta(2)}\prod_{p \mid d}\frac{p}{p+1} + O_{\epsilon}(X^{13/18+\epsilon}d^{1/4 + \epsilon}).$$
\end{lemma}

\begin{proof}
The proof for this lemma is actually simpler than that of Lemma \ref{smooth_l} because the sum over $\chi_K(n)$ only has a significant contribution when $n=1$ as opposed to whenever $n$ is a square.
\end{proof}

Turning our attention back to (\ref{to_swap}), we apply Lemma \ref{diadic_sum} by summing diadically to obtain
$$X^{1/2}\frac{\mu(\Sigma_{\Q, 2^{2i},a})}{2\zeta(2)q^2}\prod_{p \mid d}\frac{p}{p+1} + O_{\epsilon}\left(\frac{X^{13/36+ \epsilon}}{q^{43/36 - \epsilon}}\right).$$
Now we sum over $a$ and expand the definitions for $\mu(\Sigma_{\Q, 2^{2i}},a)$ and $\mu(\Sigma_{K,2^{2i}})$, which yields
\begin{align}
\sum_{a \in \{0,1,4,5\}}\mu(\Sigma_{\Q, 2^{2i},a}) \mu(\Sigma_{K,2^{2i}}) &=\left(\sum_{K_2 \in \Sigma_{\Q, 2^{2i}}}\frac{1}{\Aut{K_2/\Q_2}D_{K_2/\Q_2}}\left(\prod_{\mathfrak{p}\mid 2}\sum_{L_{\mathfrak{p}}\in \Sigma_{K,\mathfrak{p}}} \frac{1}{\Aut{L_{\mathfrak{p}}/K_{\mathfrak{p}}}D_{L_{\mathfrak{p}}/K_{\mathfrak{p}}}}\right)\right)\notag \\
&= \sum_{\substack{(L_2,K_2)\in \Sigma \\J_2(L) = 2^{2i}}}\frac{1}{\Aut{L_2,K_2}C(L_2,K_2)}. \label{sum_over_cond}
\end{align}
We will denote this by as $\mu(\Sigma_{2^{2i}})$
Bringing this back into (\ref{after_e}), we get
\begin{equation*}
    X \frac{\mu(\Sigma_{2^{2i}})}{3\zeta(2)^2}\prod_{p \mid d}\frac{1}{(p+1)^2} + O_{\epsilon}\left(\frac{X^{11/12 + \epsilon}}{q^{7/3 - \epsilon}}\right).
\end{equation*}
Summing over $q$, we get
\begin{equation*}
    X \frac{1}{3}\cdot \frac{9}{16}\left(\sum_{i = 0}^3\mu(\Sigma_{2^{2i}})\right)\prod_{p\ne 2}\left(1 - \frac{1}{p^2} - \frac{2}{p^3} + \frac{2}{p^4}\right) + O_{\epsilon}(X^{11/12 + \epsilon}).
\end{equation*}
Lastly, we note the sum $\sum \mu(\Sigma_{2^{2i}})$ is simply
$$\sum_{(L_2,K_2)}\frac{1}{\Aut{L_2,K_2}C(L_2,K_2)}.$$
Moreover, Theorem 3 of \cite{ASVW} implies,
$$\left(1 - \frac{1}{p}\right)^2\sum_{(L_p,K_p)}\frac{1}{\Aut{L_2,K_2}C(L_2,K_2)} = \left(1 - \frac{1}{p^2} - \frac{2}{p^3} + \frac{2}{p^4}\right).$$
Thus, our second sum comes to
\begin{equation}\label{easy_q}
    X \frac{3}{4}\prod_{p}\left(1 - \frac{1}{p^2} - \frac{2}{p^3} + \frac{2}{p^4}\right) + O_{\epsilon}(X^{11/12 + \epsilon}).
\end{equation}

\subsection{The Third Sum}
Our treatment of the third sum of (\ref{the_sum}) will look very similar to the second sum. We begin in the same way by summing over valid values for $q$ and get
\begin{equation}\label{third_sum}
    \sum_{\substack{q=2^i d < X^{1/4}\\ d \text{ odd,}\sqfree \\i \in \{0, 1, 2, 3\}}}\sum_{\substack{[K:\Q]=2\\ X^{1/2}/q^2 \le |D_K| < X^{1/2}}}\sum_{\substack{[L:K]=2 \\ D_{L/K} < X/|D_K| \\ J(L) = q^2}} 1.
\end{equation}
We treat the inner most sum identically by using the inclusion/exclusion principle to count only the quartic fields $L$ with $J(L) = q^2$. We also use (\ref{simple_p} -- \ref{dividing_e}) with the end result

\begin{equation}\label{after_tricks}
    X\frac{2}{3\zeta(2)}\prod_{p \mid d} \frac{1}{p(p+1)} \sum_{\substack{[K:\Q]=2\\X^{1/2}/q^2\le |D_K| < X^{1/2}\\(D_K,d)=1\\K_2 \in \Sigma_{\Q, 2^{2i}}}} \frac{\mu(\Sigma_{K,2^{2i}})}{2^{i(K)}|D_K|} \sum_{\substack{n \ge 1\\(n,2d)=1\\n \sqfree}} \frac{\chi_K(n)}{f(n)} + O_{\epsilon}\left(\frac{X^{11/12 + \epsilon}}{q^{4/3 - \epsilon}}\right).
\end{equation}

We will again sum over congruence conditions mod 8 as we did to analyze \ref{to_swap} and also use partial summation. We define
$$A(x) = \sum_{\substack{[K:\Q]=2\\0< D_K < X\\(D_K,d)=1\\D_K \equiv a \mod{8}}} \sum_{\substack{n \ge 1\\(n,2d)=1\\n \sqfree}} \frac{\chi_K(n)}{f(n)}.$$ Then, considering first the sum over real quadratic fields $K$

\begin{align*}
    \sum_{\substack{[K:\Q]=2 \\ X^{1/2}/q^2 \le D_K < X^{1/2} \\ (D_K,d) = 1\\D_K \equiv a \mod{8}}} \frac{1}{D_K}\sum_{\substack{n=1 \\ (n,2d)=1 \\ n\, \sqfree}}^{\infty} \frac{\chi_K(n)}{f(n)} &= \frac{A(t)}{t}\Big|_{X^{1/2}/q^2}^{X^{1/2}} + \int_{X^{1/2}/q^2}^{X^{1/2}} \frac{A(t)}{t^2}dt \\
    &= \frac{\mu(\Sigma_{\Q, 2^{2i},a})\log(q^2)}{3\zeta(2)}\prod_{p \mid d} \frac{p}{p+1} + O_{\epsilon}\left(\frac{q^{29/36+\epsilon}}{X^{5/36 - \epsilon}}\right).
\end{align*}

As before, we extend the sum over the different congruence conditions and imaginary quadratic fields and bring this back into (\ref{after_tricks}) to get

\begin{equation}\label{last_piece}
    X \frac{1}{3\zeta(2)^2}\sum_{\substack{q=2^i d < X^{1/4}\\ d \text{ odd,}\sqfree \\i \in \{0, 1, 2, 3\}}}\mu(\Sigma_{2^{2i}})\log(q^2)\prod_{p \mid d} \frac{1}{(p+1)^2} + O_{\epsilon}(X^{11/12 + \epsilon}).
\end{equation}
If we replace $q$ with $2^i d$, we can split $\log(q^2)$ and consider the sum in two parts. First, we have
\begin{equation*}
    \sum_{\substack{q=2^i d < X^{1/4}\\ d \text{ odd,}\sqfree \\i \in \{0, 1, 2, 3\}}}\mu(\Sigma_{2^{2i}})\log(d^2)\prod_{p \mid d} \frac{1}{(p+1)^2} = 5\sum_{d \text{ odd,}\sqfree}\log d\prod_{p \mid d}\frac{1}{(p+1)^2} + O_{\epsilon}(X^{-1/4 + \epsilon}).
\end{equation*}
We can rewrite this as a sum over odd primes $p$ and get
\begin{align*}
    5\sum_{d \text{ odd,}\sqfree}\log d \prod_{p\mid d}\frac{1}{(p+1)^2} &= 5\sum_{p\ne 2}\frac{\log p}{(p+1)^2}\prod_{r \ne p, 2}\left(1 + \frac{1}{(r+1)^2}\right)\\
    &= \frac{9}{4}\prod_{p}\left(1 + \frac{1}{(p+1)^2}\right)\left(-\frac{\log 2}{5} + 2\sum_{p} \frac{\log p}{p^2 + 2p + 2}\right).
\end{align*}
The second sum is
\begin{align*}
    2 \log 2\sum_{\substack{q=2^i d < X^{1/4}\\ d \text{ odd,}\sqfree \\i \in \{0, 1, 2, 3\}}} i\cdot \mu(\Sigma_{2^{2i}}) \prod_{p \mid d}\frac{1}{(p+1)^2} &= \frac{9\log 2}{5}\prod_{p}\left(1 + \frac{1}{(p+1)^2}\right)\left(\sum_{i = 1}^3 i \cdot\mu(\Sigma_{2^{2i}}) \right) + O(X^{-1/4}).
\end{align*}

Using Tables 1-4, we compute $\sum i\cdot \mu(\Sigma_{2^{2i}}) = 11/16$ and put everything back into (\ref{last_piece}).
\begin{equation}\label{ready_to_go}
    X\frac{3}{4}\left(\frac{7 \log 2}{20} + 2\sum_p \frac{\log p}{p^2 + 2p + 2}\right)\prod_{p}\left(1 - \frac{1}{p^2} - \frac{2}{p^3} + \frac{2}{p^4}\right) + O_{\epsilon}(X^{11/12 + \epsilon}).
\end{equation}

The main theorem now follows from putting (\ref{first_sum}), (\ref{easy_q}), and (\ref{ready_to_go}) in the sum (\ref{the_sum}).

\subsection{A Modified Argument for Theorem \ref{local_specs}}
To prove Theorem \ref{local_specs}, one can follow the same outline as the main theorem by tackling each sum in (\ref{the_sum}) separately. Let $\Sigma = (\Sigma_v)_v$ be an acceptable collection of local specifications and let $m$ be the product of primes $p$ such that $\Sigma_p$ doesn't contain every pair $(L_p, K_p).$ If we again restrict $\Sigma$ to $\Sigma_K$, then Theorem \ref{local_asymptotics} gives
\begin{equation*}
    \sum_{\substack{[L:K]=2\\ L \in \mathcal{K}(\Sigma_K) \\ D_{L/K} \le X/|D_K|}} = X \frac{L(1,K/\Q)}{|D_K|L(2,K/\Q)} \cdot \mu(\Sigma_{K,\infty}) \prod_{p \mid m} \left(\prod_{\mathfrak{p}\mid p} \frac{\mathrm{N}\mathfrak{p}}{1 + \mathrm{N}\mathfrak{p}}\right) \cdot \mu(\Sigma_{K, p}) + O_{\epsilon}(X^{2/3 + \epsilon}|D_K|^{-1/2 + \epsilon}m^{1/3 + \epsilon}),
\end{equation*}
where
$$\mu(\Sigma_{K, \infty}) = \prod_{\sigma_v \mid \infty}\sum_{L_v \in \Sigma_{K, v}} \frac{1}{\Aut{L_v/K_{v}}}, \text{ and } \mu(\Sigma_{K,p})=\prod_{\mathfrak{p}\mid p  }\sum_{L_{\mathfrak{p}} \in \Sigma_{K, \mathfrak{p}}} \frac{1}{\Aut{L_\mathfrak{p}/K_\mathfrak{p}}D_{L_\mathfrak{p}/K_\mathfrak{p}}}.$$
The $m^{1/3+\epsilon}$ is a worst case estimate in which every prime $p$ dividing $m$ were either split or inert in $K$ and $\mathfrak{p} \mid \mathfrak{u}_1$ for every $\mathfrak{p} \mid p$ in our application of Theorem \ref{local_asymptotics}.

At this point, we use (\ref{simple_p}) and note that we can pull these factors outside the sum over $K$ as well as $\mu(\Sigma_{K,\infty})$ and $\mu(\Sigma_{K,p})$ and have
\begin{equation*}
    X \frac{\mu(\Sigma_{K,\infty})}{\zeta(2)}\prod_{p\mid m}\left(\left(1 + \frac{1}{p}\right)^{-1}\mu(\Sigma_{K,p})\right) \sum_{\substack{[K:\Q]=2 \\K \in \mathcal{K}(\Sigma_{\Q}) \\ |D_K| \le X^{1/2}}} \frac{1}{|D_K|} \prod_{p \nmid m}\left(\frac{1 - \chi_K(p)/p^2}{1 - \chi_K(p)/p}\right) + O_{\epsilon}(X^{11/12 + \epsilon}m^{1/3 + \epsilon}).
\end{equation*}

To finish the first sum, only slight modifications to Lemmas \ref{theorem_2} and \ref{smooth_l} are needed. A modification worth mentioning is the one we need to define $c_{\Sigma}$ in the theorem statement. In the modified Lemma \ref{theorem_2} we would define

$$S_{\Sigma}^{\pm}(X) = \sum_{\substack{[K:\Q]=2 \\ K \in \mathcal{K}(\Sigma_{\Q}) \\0< \pm D_K< X}} \frac{1}{|D_K|}\cdot \frac{L(1,K/\Q)}{L(2,K/\Q)}.$$
Then,
$$E_{\Sigma}^{\pm}(X) = S_{\Sigma}^{\pm}(X) - \frac{1}{2}X \prod_p \left(\left(1 - \frac{1}{p}\right)^2\sum_{(L_{p},K_{p})\in \Sigma_{p}} \frac{1}{\Aut{L_{p},K_{p}}C(L_p,K_p)}\right).$$
With $c_{\Sigma}^{\pm}$ defined analogously to $c^{\pm}$, we have
\begin{equation}\label{local_c}
    c_{\Sigma} = c_{\Sigma}^+ + \frac{1}{2}c_{\Sigma}^-.
\end{equation}

Apart from this, we should also note that given an acceptable collection of local specifications $\Sigma$ that is not complete at some prime $p$, it is not clear that doubling the first sum yields that correct result as the second iteration of the sum is meant to count pairs $(L, K)$ for which $D_{\phi(K)} < X^{1/2}$ and $(\phi(L),\phi(K))$ may not be in $\mathcal{L}(\Sigma).$ But, the discussion in section 9.2 of \cite{ASVW} shows that $\mu(\Sigma_p) = \mu(\phi(\Sigma_p))$ if we give $\phi(\Sigma_p)$ the meaning you might expect. 

For the remaining two sums, the same sorts of modifications can be made without much change to the original arguments. Moreover, as each sum is only computed once, we don't need to worry about the equivalence of $\mu(\Sigma_p)$ and $\mu(\phi(\Sigma_p)).$ However, we will note that the constant $\frac{7 \log 2}{20}$ in the secondary term of Theorem \ref{main_thm} is
$$2 \log 2\frac{2\mu(\Sigma_{2^4}) + 3\mu(\Sigma_{2^6})}{\mu(\Sigma_2)}$$
and, thus, agrees with Theorem \ref{local_specs} when $\Sigma_2$ contains all pairs $(L_2, K_2)$.

\end{document}